\documentclass[final,leqno,onefignum,onetabnum]{siamltex1213}
\pdfoutput=1
\title{
Robustness analysis of spatiotemporal models in the presence of extrinsic fluctuations
}

\author{Andreas Hellander \thanks{Division of Scientific Computing,
			Department of Information Technology,
			Uppsala University, P. O. Box 337, SE-75105 Uppsala, Sweden. (\email{ {andreas.hellander@it.uu.se}, {perl@it.uu.se}}) } \and Jan Klosa  \footnotemark[1] \and Per L\"{o}tstedt \footnotemark[1]  \and Shev MacNamara \thanks{School of Mathematics and Statistics, University of New South Wales, Sydney, Australia. (\email{ {s.macnamara@unsw.edu.au})} }
}

\usepackage{amsmath}
\usepackage{subfigure}
\usepackage{url}
\usepackage{doi}
\usepackage[numbers]{natbib}
\usepackage[margin=3.5cm]{geometry}
\usepackage{amssymb}
\usepackage{graphicx}

\def\veps{\varepsilon}
\def\fatf{\mathbf{f}}
\def\fatg{\mathbf{g}}
\def\fath{\mathbf{h}}
\def\fatn{\mathbf{n}}
\def\fats{\mathbf{s}}
\def\fatu{\mathbf{u}}
\def\fatx{\mathbf{x}}

\def\fatvarsigma{\boldsymbol{\varsigma}}

\def\fatsigma{\boldsymbol{\sigma}}

\def\fatrho{\boldsymbol{\rho}}
\def\hdfr{\delta\hat{\fatrho}}
\def\dfr{\delta\fatrho}
\def\hdfr{\delta\hat{\fatrho}}
\def\dr{\delta\rho}
\def\hdr{\delta\hat{\rho}}

\def\hdfu{\delta\hat{\fatu}}
\def\tdfu{\delta\tilde{\fatu}}

\def\du{\delta u}
\def\hdu{\delta\hat{u}}
\def\tdu{\delta\tilde{u}}
\def\hG{\hat{G}}
\def\tG{\tilde{G}}
\def\hH{\hat{H}}
\def\ps{\partial_s}
\def\pt{\partial_t}
\def\px{\partial_x}
\def\px{\partial_x}

\def\ordo{\mathcal{O}}

\newcommand{\ud}{\mbox{d}}

\renewcommand{\figurename}{Fig.}
\newcommand{\captionfonts}{\small}

\makeatletter  
\long\def\@makecaption#1#2{%
  \vskip\abovecaptionskip
  \sbox\@tempboxa{{\captionfonts #1: #2}}%
  \ifdim \wd\@tempboxa >\hsize
    {\captionfonts #1: #2\par}
  \else
    \hbox to\hsize{\hfil\box\@tempboxa\hfil}

  \fi
  \vskip\belowcaptionskip}
\makeatother   

\begin{document}

\maketitle


\begin{abstract}
We analyze the governing partial differential equations of a model of pole-to-pole oscillations of the \textit{MinD} protein in a bacterial cell.
The sensitivity to extrinsic noise in the parameters of the model is explored.
Our analysis shows that overall, the oscillations are robust to extrinsic perturbations in the sense that small 
perturbations in reaction coefficients result in small differences in the frequency and in the amplitude.
However, a combination of analysis and simulation also reveals that the oscillations are more sensitive to some extrinsic time-scales than to others.

\end{abstract}



\section{Introduction}
Mathematical models are now essential to the way biological scientists understand single cells~\citep{GunawardenaPatheitc,Kirschner50YearsJacobMonod2011,PhilipsPhysBiologyCellBook2012}.
Chemical reactions and transport of chemical species are often described by deterministic models, for example, by partial differential equations (PDEs) and the Law of Mass Action.
However, noise plays a fundamental role in many cellular processes~\citep{BlaKae03,FedFon02,ShaSwa08} such as switching between stable modes of gene expression~\citep{FanElf06,McAArk97,VandenFabio2006,TiaBur06}.
For such processes, a discrete and stochastic modelling framework is more appropriate than a deterministic continuum model, especially when a single cell contains only a small number of molecules of a particular chemical species~\citep{BlaKae03,FedFon02,Kam01}.
Such a framework is provided by the chemical master equation \citep{Gil92,Kam01}, which is increasingly applied in systems biology.

When formulating a stochastic model of a process, we may distinguish between external or \textit{extrinsic noise} that is independent of the system being modeled, and internal or  \textit{intrinsic noise} that is inherently part of the system itself.
Van Kampen discusses this issue in his classic text.
He emphasizes the importance of making this distinction at a conceptual level during the process of model formulation~\citep[Chapter IX.5, Chapter XVII.7]{Kam01}.
Biologists also employ the terminology of intrinsic and extrinsic noise when describing stochastic phenomena in relation to models of gene expression, although identifying and measuring intrinsic and extrinsic contributions to dynamic systems can be challenging~\citep{EloLev02,HilfingerPaulsson2011,SwaElo02}.
A common interpretation is that intrinsic noise arises from the inherently discrete nature of a collision theory of chemical reactions, in which there is randomness associated with the chance collisions of molecules, whereas extrinsic noise arises from all of the other processes that we do not explicitly include in the mechanistic steps of our mathematical model but which we do believe exert influence.
The stage of the cell-cycle, ambient temperature, a dynamic microenvironment, or the number of ribosomes in a cell, all effect cellular processes but they are usually not explicitly included in models; instead their effects may be regarded as extrinsic noise.
For example, temperature effects chemical reaction rates and also biological oscillations \citep{Gould13,TouJerRut:2006,WaAnDiCu}.


We are interested in a mathematical model of MinD oscillations in bacteria \citep{FanElf06,KruHowMar:2007,RenRuo02,WaAnDiCu}
and the robustness of the model to extrinsic spatial and temporal fluctuations in the coefficients.
The model is a system of nonlinear PDEs with diffusion for the mean values of the concentrations of the species.
If the copy number of the molecular species is large then the relative intrinsic fluctuations are small and a deterministic PDE system without intrinsic noise is a good approximation.
This is often the case for the MinD oscillations \citep{Kerr06,WaAnDiCu} but not in all situations \citep{FanElf06}. The question of robustness is certainly important specifically in the context of models of MinD oscillations \citep{Halatek12,KruHowMar:2007},
but it is also important more generally in the field of uncertainty quantification and in systems biology, where parameters are often poorly characterised \citep{WeisseMiddletonHuisinga10}.
Moreover, oscillations in biology have a rich literature, in which robustness of oscillations is an important theme, e.g.
\citep{Mur02,WinfreeGeometryOfBiologicaltime}.
With our approach combining analysis with simulations we find that overall, the MinD model is robust to fluctuations in the coefficients,  in the sense that small
fluctuations in the coefficients lead to only small changes in the period or to small changes in the amplitude of the oscillations.
However, our results also reveal that the oscillations are more sensitive to some timescales of the extrinsic fluctuations than to others.

The outline of the paper is as follows. The PDE model of the Min oscillations is found in Section~\ref{sec:Minmodel}. The model is linearized, expanded in a cosine series, and in a small
parameter $\veps$ in Section~\ref{sec:macromodel}. The extrinsic perturbations of the parameters of the model are scaled by $\veps$ and
the influence of the perturbations on the frequency and the amplitude of the oscillations is analyzed.
Section~\ref{sec:extnoise} is a brief review of the properties of an Ornstein-Uhlenbeck process for the colored temporal noise and how the spatial noise is generated. The autocorrelations for
the changes in frequency are derived in Section~\ref{sec:randpert} assuming that the perturbations are as in Section~\ref{sec:extnoise}.
Comparison is made in Section~\ref{sec:comparisons} between the solutions of the nonlinear model and the linearized model used in the analysis.
Some conclusions are drawn in the final Section~\ref{sec:disc}.

\begin{figure}[h]
\centering
\subfigure{\includegraphics[height=8cm,width=0.45\linewidth]{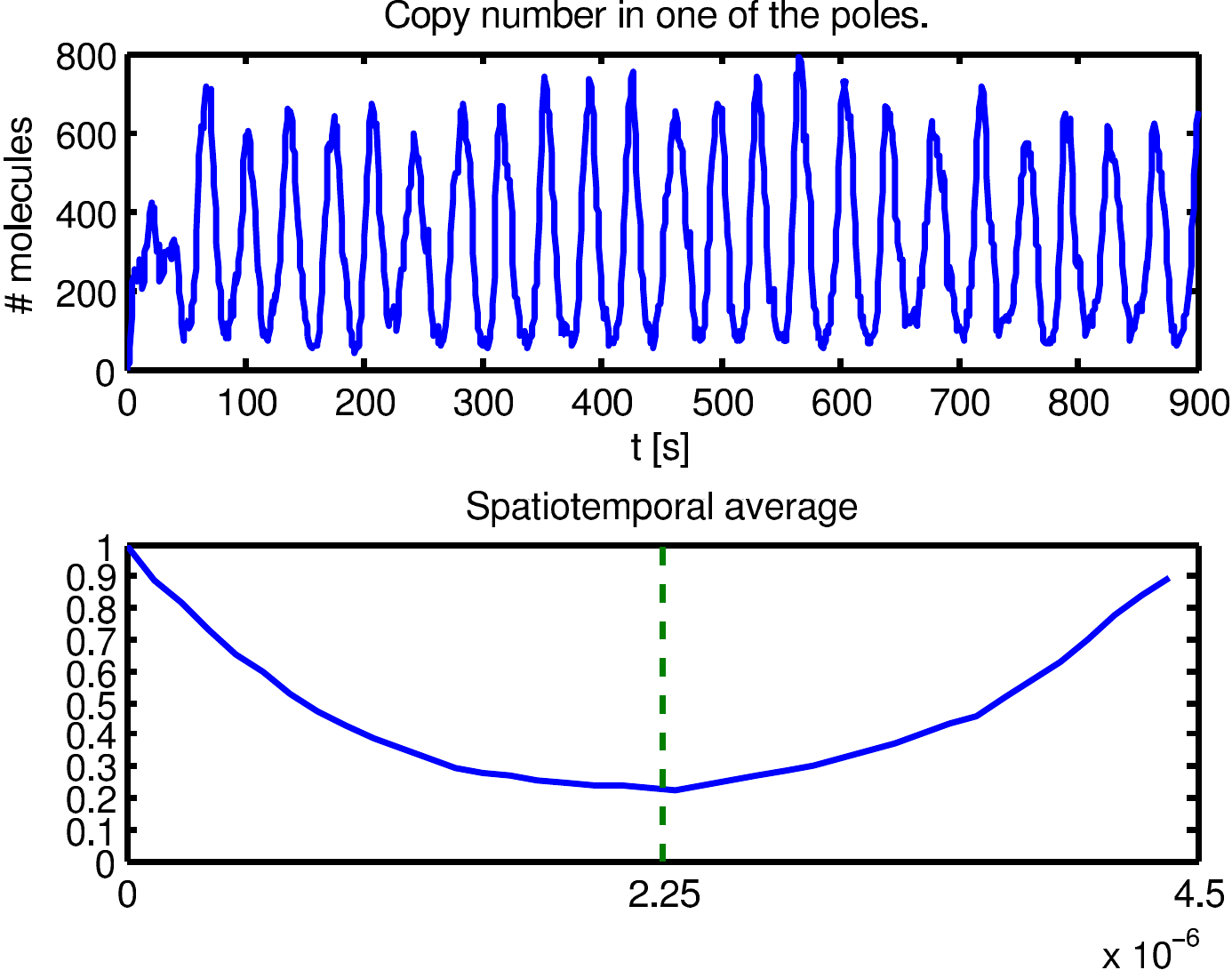}}
\subfigure{\includegraphics[height=8cm,width=0.45\linewidth]{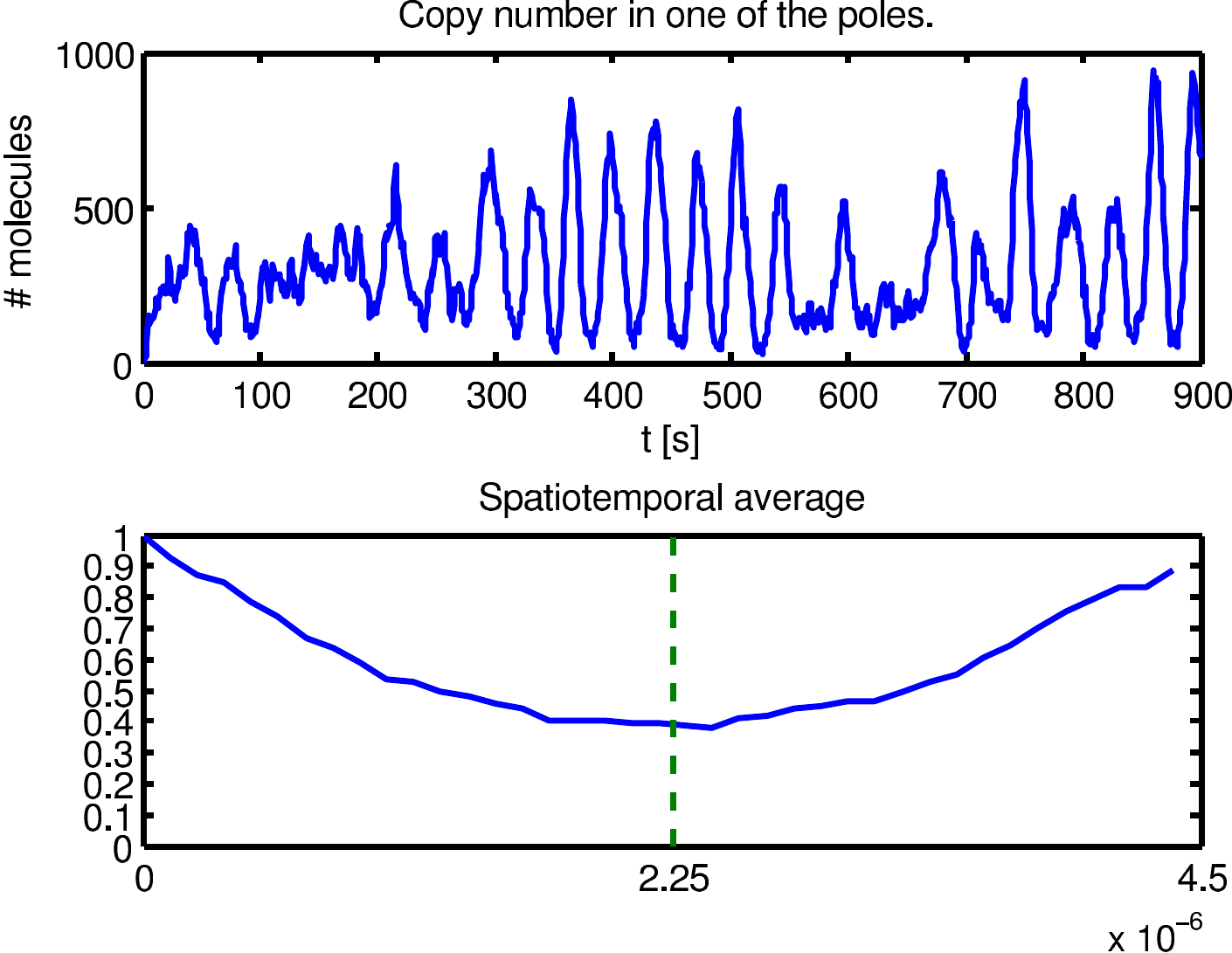}}
\caption{Stochastic simulations in three space dimensions of the oscillations in MinD proteins at one pole (top) and the time-averaged concentration profile along an {\it E. coli} bacterium (bottom)
with unperturbed parameters (left) and temporal and spatial perturbations of the parameters (right).
}
\label{fig:URDMETimeSeries}
\end{figure}

 \begin{figure}[htbp]
\centering
\subfigure[]{\includegraphics[height=6cm,width=\linewidth]{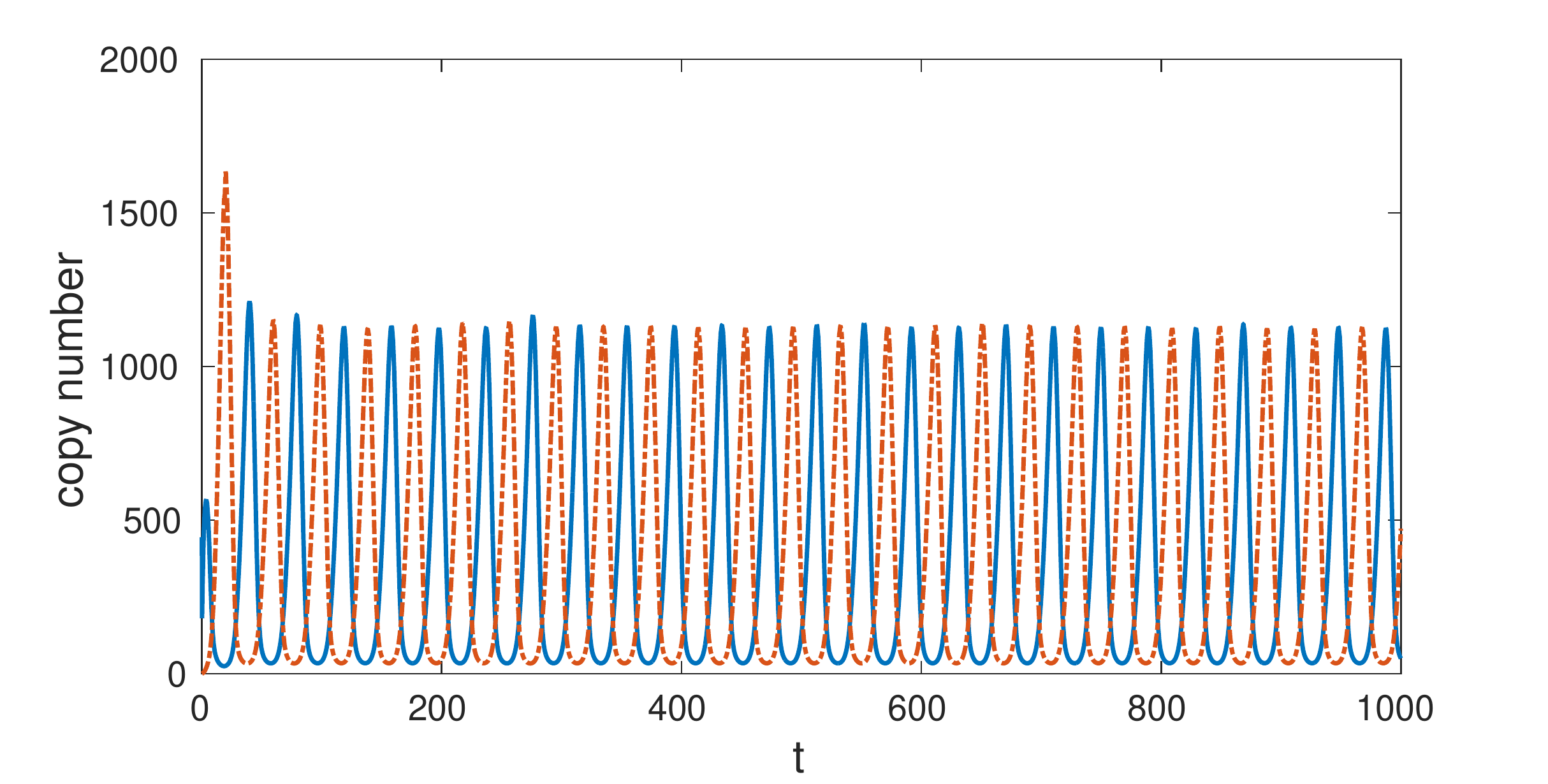}}
\subfigure[]{\includegraphics[width=\linewidth]{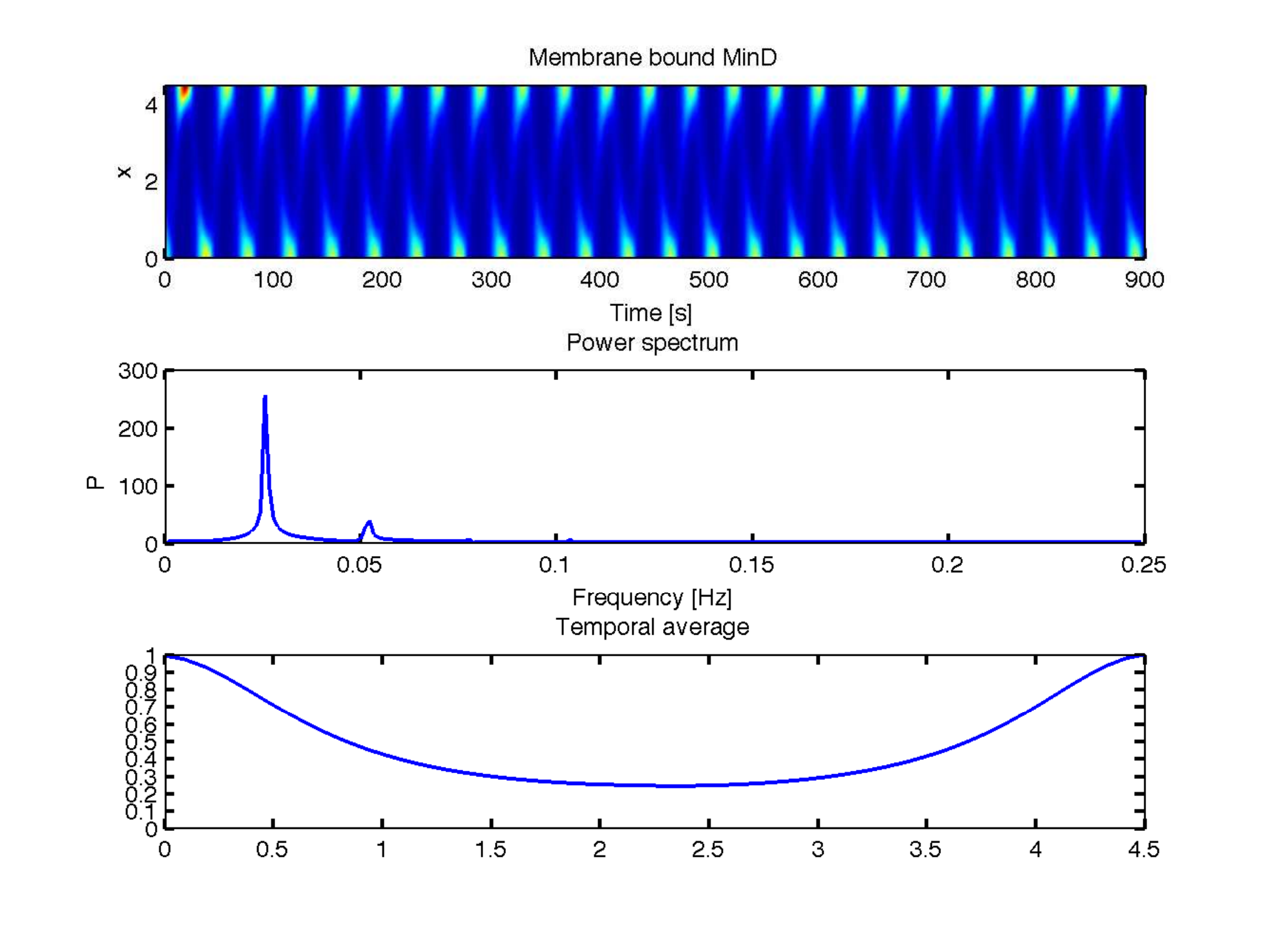}}
\caption{Solution of the deterministic PDE model, in one space dimension and time, of MinD oscillations \eqref{eq:Huang_model}.
(a) Temporal variation of $\rho_d$ shown at the left boundary at $x=0$ (solid blue) and at the right boundary $x=4.5$ (dashed red).
(b) Top: Kymographs of the dynamics of the concentration show regular oscillations (high concentration: red, low concentration: blue).
Most time is spent at the poles, with a relatively fast transition from one pole to the other.
Middle: Power spectrum confirms strong periodic components of the solution.
The maximal peak corresponds to a period of about 40 s.
Bottom: The time-averaged concentration profile is at a minimum near the middle of the cell.
}
\label{fig:mincdepde}
\end{figure}
\section{MinD proteins oscillate in a single cell}\label{sec:Minmodel}

 Experimental observations of a single bacteria cell reveal that MinD proteins oscillate from one pole of the cell to the other, 
with a period of about one minute \citep{FanElf06,Kruse:2002,KruHowMar:2007,MeaKru:2005,RenRuo02,WaAnDiCu}.
 During these oscillations, Min proteins spend most of the time at the poles of the cell and much less time at the middle of the cell, so that a time-averaged profile shows MinD concentration lowest in the middle of the cell and highest at the poles of the cell.
 These oscillations in space and time are associated with correct functioning of cell division.
 The time-averaged MinD  concentration profile can be thought of as a potential function that repels  key cellular machinery (such as Ftsz proteins and assembly of the Z-ring) from the poles of the cell, and instead pushes the machinery to the middle of the cell, where MinD concentration is lowest.
 This allows the cell to correctly locate and divide at approximately the middle, which is important for producing two equal-sized daughter cells.
  Disruptions of these oscillations are associated with cells that divide unevenly or that exhibit other problematic phenotypes \citep{FanElf06} so robustness is an important issue \citep{Halatek12, SteGilDoy2004, StellingDoyleRobustness2004}.

The Min system is simulated stochastically in three dimensions (3D) with a mesoscopic model and 
Gillespie's SSA \citep{gillespie} implemented in \citep{URDME} in 
Figure~\ref{fig:URDMETimeSeries}. 
Parameters are not perturbed in the left column, so the noise is intrinsic there. In the right column,
one parameter, $\sigma_{dD}$ in \eqref{eq:Huang_model}, is perturbed in space and time and we have both intrinsic and extrinsic noise. 
The MinD oscillations are affected by the extrinsic noise in the upper right panel but the average concentration profile is less sensitive.

A system of reaction-diffusion PDEs is a popular macroscopic model for the oscillations of the Min protein \citep{Huang:2003,KruHowMar:2007}.
It includes five species with concentrations that vary in space and time: three species in the cytosol of the cell, and two species that are membrane-bound.
Let the concentrations of MinD:ADP, MinD:ATP, and MinE in the cytosol be
$\rho_{DD}(\fatx, t), \rho_{DT}(\fatx, t),$ and $\rho_E(\fatx, t)$.
Let $\rho_{d}(\fatx, t)$ and $\rho_{de}(\fatx, t)$ be the concentrations of MinD:ATP and MinE:MinD:ATP, which are complexes on the membrane.
The volume of the domain (the cytosol in a single cell) is denoted by $\Omega$ with the boundary (the membrane of the cell) $\partial\Omega$  and an outward normal $\fatn$.

The equations for the concentrations of the species in the model of Huang \textit{et al.} \citep{Huang:2003} are
\begin{equation}
\begin{array}{lll}
  \partial_t \rho_{DD}&=\sigma_{de}\rho_{de}-\sigma_{DT}\rho_{DD}+\gamma_D\Delta \rho_{DD},\\
  \partial_t \rho_{DT}&=\sigma_{DT}\rho_{DD}-(\sigma_{D}+\sigma_{dD}(\rho_d+\rho_{de}))\rho_{DT}+\gamma_D\Delta \rho_{DT},\\
  \partial_t \rho_{E}&=\sigma_{de}\rho_{de}-\sigma_{E}\rho_{d}\rho_{E}+\gamma_E\Delta \rho_{E},\\
  \partial_t \rho_{d}&=(\sigma_{D}+\sigma_{dD}(\rho_{d}+\rho_{de}))\rho_{DT}-\sigma_{E}\rho_{d}\rho_{E},\\
  \partial_t \rho_{de}&=\sigma_{E}\rho_{d}\rho_{E}-\sigma_{de}\rho_{de}.
\end{array}
\label{eq:Huang_model}
\end{equation}
The time derivative is denoted by $\partial_t$ and the diffusion operator by $\Delta$.
Reactions involving $\rho_d$ and $\rho_{de}$ take place only on the cell membrane.
The boundary conditions for the species in the cytosol are reflective at $\partial\Omega$, i.e. $\fatn\cdot\nabla \rho=0$.

Other models of the Min system are reviewed in \citep{KruHowMar:2007}.
The model of Fange and Elf \citep{FanElf06} has diffusion also on the membrane and the term $\sigma_{dD}\rho_{de}\rho_{DT}$ is missing in the second and fourth equations in \eqref{eq:Huang_model}.
The change of MinD from ADP to ATP form is ignored and there is an upper bound on the number of membrane binding sites in the model of Meacci and Kruse \citep{MeaKru:2005}.

The geometry of the cell $\Omega$ is modelled as cylindrical, with spherical caps at both ends.
The cell radius is $0.5 \mu m$, the cylindrical part is $3.5 \mu m$, and the volume $V$ is $3.2725 \mu m^3$.
The typical reaction parameters in \eqref{eq:Huang_model} are
\begin{equation}
\begin{array}{lll}
  \sigma_{de}=0.7 s^{-1}, \sigma_{DT}=1 s^{-1}, \;\;\;\; \sigma_D=0.025 \mu m s^{-1}, \\
  \sigma_{dD}=6.8\cdot 10^5 M^{-1} s^{-1}, \;\;\;\; \sigma_{E}=5.60\cdot 10^7 M^{-1} s^{-1}.
\end{array}
\label{eq:Huang_param}
\end{equation}
The diffusion coefficients are $\gamma_D=\gamma_E=2.5 \mu m^2 s^{-1}$. The cell length is assumed to be constant although 
the length is varying during the cell cycle and has an influence on the oscillations \citep{FiFrMeLuChKr2010}.

Figure~\ref{fig:mincdepde} shows the deterministic solution of the system \eqref{eq:Huang_model} for the parameters \eqref{eq:Huang_param} in one dimension (1D).
As can be seen, periodic oscillations of MinD from pole to pole (top) works to establish a relative temporal average concentration profile in which MinD has a
higher concentration in the regions near the polar caps and a minimum in the middle of the cell (bottom).
For these values of the parameters, the power spectrum (middle pane) has its main peak at approximately 0.025 Hz, corresponding to a period of approximately 40 seconds.
The solution is similar to the stochastic realization in Figure~\ref{fig:URDMETimeSeries} with unperturbed parameters.
The steady state solution $\fatrho_\infty$ of \eqref{eq:Huang_model} agrees very well 
with the average values in space and time of the
stochastic simulations. The time period $T$ of the oscillations is about 40 seconds in both
the deterministic equations and the stochastic simulations.

\section{Analysis of the macroscopic model}\label{sec:macromodel}

The deterministic PDE model switches between a spatially homogeneous equilibrium and unstable periodic oscillations, 
via stable oscillations when parameters $\fatsigma$ are varied.
The oscillatory behavior can be compatible with the suppression of Z-ring formation at the bacterium's poles only if 
the oscillations are reliable enough in space and of large enough amplitude.
We are interested in oscillatory solutions of \eqref{eq:Huang_model} and where in the parameter space they appear.
Small perturbations around a steady state solution $\fatrho_\infty$ (or a fixed point) are introduced.
The small perturbations satisfy linearized equations with a constant system matrix.
The eigenvalues of this matrix tell us where the perturbations are stable, unstable, or oscillate.
The coefficients $\fatsigma$ are perturbed in space and time
about a constant mean value. In this way, the uncertainty in the parameters is introduced. The amplitude and 
the frequency of the oscillations in the MinD system
are changed by the perturbations which are assumed to be small such that
linearization is possible.

\subsection{Invariants in the deterministic model}

Since
\[
   \int_{\Omega}\Delta \rho\,\ud\Omega=\int_{\partial\Omega} \fatn\cdot\nabla\rho\, \ud S=0,
\]
it follows from \eqref{eq:Huang_model} that the
total number of MinD and MinE molecules, $N_D$ and $N_E$, defined by
\begin{equation}
\begin{array}{lll}
  N_D&=\int_{\Omega}\rho_{DD}+\rho_{DT}+\rho_{d}+\rho_{de}\, \ud\Omega,\quad
  N_E&=\int_{\Omega}\rho_{E}+\rho_{de}\, \ud\Omega,
\end{array}
\label{eq:Huang_total}
\end{equation}
are constant and $\pt N_D=\pt N_E=0$. The total number of molecules in our examples are
$N_D=4500$ and $N_E=1575$.

The conclusion from \eqref{eq:Huang_total} for a constant steady state solution
\begin{equation}
 \fatrho_{\infty}^T=(\rho_{DD\infty}, \rho_{DT\infty}, \rho_{E\infty}, \rho_{d\infty}, \rho_{de\infty}),
\label{eq:Huang_stst}
\end{equation}
is that the quantities $\rho_{Dtot}$ and $\rho_{Etot}$ in
\begin{equation}
  \rho_{Dtot}=N_D/V=\rho_{DD\infty}+\rho_{DT\infty}+\rho_{d\infty}+\rho_{de\infty}, \quad
  \rho_{Etot}=N_E/V=\rho_{E\infty}+\rho_{de\infty},
\label{eq:Huang_cons}
\end{equation}
are conserved in all solutions. Then $\rho_{d\infty}$ and $\rho_{de\infty}$ can be eliminated from
the stationary equation of \eqref{eq:Huang_model} using \eqref{eq:Huang_cons} yielding three nonlinear equations
for $\rho_{DD\infty}, \rho_{DT\infty},$ and $\rho_{E\infty}$. The constants $\rho_{Dtot}$ and
$\rho_{Etot}$ are $1375$ and $481$. Only one fixed point has been found in the neighborhood
of the $\fatsigma$-values in \eqref{eq:Huang_param} and it depends smoothly on the parameters.

\subsection{Model with variable parameters}

In order to investigate the influence of a variation in the $\fatsigma$-parameters, a 1D
simplification of the model in \eqref{eq:Huang_model} is
introduced in the interval $[0, L]$ with $L=4.5 \mu m$ and
\begin{equation}
   \pt\fatrho=\fatf(\fatrho)+\kappa(x, t)\fatg(\fatrho)+\gamma D\partial^2_x\fatrho,\quad \px \fatrho=0\; {\rm at}\;x=0, L .
\label{eq:Huang_1D}
\end{equation}
Here $\fatf$ and $\fatg$ contain the reaction terms and $D$ is diagonal with $D_{jj}=1, j=1,2,3,$ and $D_{jj}=0, j=4,5$ in the MinD model,
and $\px$ denotes $\partial/\partial x$. A 1D model is found to be sufficient to study Min oscillations in \citep{KruHowMar:2007}. 
See also Figures~\ref{fig:URDMETimeSeries} and \ref{fig:mincdepde}.
The parameters \eqref{eq:Huang_param} in \eqref{eq:Huang_model} are constant in $\fatf$ and are multiplied by the same factor
$\kappa$ varying in space and time in $\fatg$. The assumption is that the perturbed parameters appear linearly in the right hand
side of \eqref{eq:Huang_1D}. The factor is assumed to have the expansion
\begin{equation}
\begin{array}{lll}
   \kappa(x, t)&=(1+\veps \kappa_{t}(t)+\ordo(\veps^2))(1+\veps\kappa_x(x)+\ordo(\veps^2))\\
               &=1+\veps \kappa_{t}(t)+\veps\kappa_x(x)+\ordo(\veps^2)
\end{array}
\label{eq:kappadef}
\end{equation}
in a small parameter $\veps$. The perturbations are such that
\begin{equation}
   \lim_{t\rightarrow \infty}\frac{1}{t}\int_0^T \kappa_t(s)\, \ud s=0, \quad \frac{1}{L}\int_0^L \kappa_x(x)\, \ud x=0.
\label{eq:enoise}
\end{equation}
Thus, the mean values of the $\fatsigma$-parameters in space and time are not changed.
The unperturbed constant steady state $\fatrho_{\infty}$ with $\kappa=1$ satisfies
\begin{equation}
   \fatf(\fatrho_{\infty})+\fatg(\fatrho_{\infty})=0.
\label{eq:ststate}
\end{equation}
A perturbation of $\fatrho_{\infty}$ is denoted by
\begin{equation}
  \dfr(x, t)^T=(\delta\rho_{DD}(x, t), \delta\rho_{DT}(x, t),
  \delta\rho_{E}(x, t), \delta\rho_{d}(x, t), \delta\rho_{de}(x, t)).
\label{eq:Huang_pert}
\end{equation}
Insert $\fatrho_{\infty}+\dfr$ into \eqref{eq:Huang_1D} and linearize the system of
equations. Terms of $\ordo(\|\dfr\|^2)$ are ignored and the Jacobians of $\fatf$ and $\fatg$ at $\fatrho_\infty$
are denoted by $F=\partial\fatf/\partial\fatrho$ and $G=\partial\fatg/\partial\fatrho$. Then $\dfr(x, t)$ satisfies
\begin{equation}
\begin{array}{lll}
   \pt\dfr&=F\dfr +\kappa(x, t)G \dfr+\gamma D\px^2\dfr\\
                                &=J \dfr+\gamma D\px^2\dfr
                                 +\veps\kappa_{t}(t)G \dfr +\veps\kappa_x(x)G \dfr +\ordo(\veps^2).
\end{array}
\label{eq:Huang_eps}
\end{equation}
where $J=F+G$. The constant Jacobian matrix $J$ depends on the steady state solution $\fatrho_{\infty}$ and the unperturbed reaction coefficients.
The expansion of $\dfr$ in the small parameter is
\begin{equation}
   \dfr=\dfr_0+\veps\dfr_1+\veps^2\dfr_2+\ordo(\veps^3).
\label{eq:drhoeps}
\end{equation}
This expansion will be inserted into \eqref{eq:Huang_eps} to derive equations for $\dfr_0$ and $\dfr_1$ but
first the stability of the lowest order term in the expansion is investigated.

\subsection{Stability analysis of the lowest perturbation mode}\label{sec:stabanal}

An equation for the unperturbed solution $\dfr_0(x, t)$ is obtained by letting $\veps=0$ in \eqref{eq:Huang_eps}
\begin{equation}
   \pt\dfr_0=J \dfr_0+\gamma D\px^2\dfr_0,
\label{eq:Huang_dfr0}
\end{equation}
satisfying the constraints obtained from \eqref{eq:Huang_total}
\begin{equation}
   \delta\rho_{DD}+\delta\rho_{DT}+\delta\rho_{d}+\delta\rho_{de}=0,\; \delta\rho_{E}+\delta\rho_{de}=0.
\label{eq:Huang_constr}
\end{equation}

The stability of the constant steady state is first investigated by letting $\dfr_0$ be constant
in space in \eqref{eq:Huang_dfr0}, $\dfr_0=\dfr_0(t)$. Then the equation for $\dfr_0$ is
\begin{equation}
   \pt\dfr_0=J \dfr_0,
\label{eq:Huang_dfr00}
\end{equation}
with the solution
\begin{equation}
   \dfr_0(t)=\exp(Jt) \dfr_0(0).
\label{eq:Huang_sol0}
\end{equation}
The eigenvalues $\lambda_j(J)$ of $J$ determine the stability properties of the solution. If
$\max\Re\lambda_j< 0$ then a spatially constant perturbation will vanish but
suppose that $\max\Re\lambda_j>0$ and $\dfr_0(0)\ne 0$. Then there are growing perturbations
violating the assumption of small perturbations. Furthermore, in order to satisfy \eqref{eq:Huang_constr} at least
one component must approach $-\infty$ breaking the non-negativity constraint on the concentrations.
Therefore, we let the constant steady state be unperturbed initially with $\dfr_0(0)= 0$ and $\dfr_0(t)= 0$.


Another perturbation mode satisfying the boundary conditions in \eqref{eq:Huang_1D} is
\begin{equation}
   \dfr_0(x, t)=\dfr_0(t)\cos(\pi x/L).
\label{eq:drhocos}
\end{equation}
The solution to \eqref{eq:Huang_dfr0} with this ansatz is
\begin{equation}
   \dfr_0(x, t)=\exp(H_1t)\cos(\pi x/L)\dfr_0(0),\quad  H_1=J-(\gamma\pi^2/L^2) D.
\label{eq:drhosol}
\end{equation}
The stability of this perturbation is determined by the eigenvalues $\lambda_j(H_1), j=1,\ldots,5$. In
the neighborhood of the
$\fatsigma$-values in \eqref{eq:Huang_param}, there is an oscillatory mode with $\Re\lambda_{j}(H_1)=0,\, j=1,2,$ and $\Re\lambda_{j}(H_1)<0, j=3,4,5,$
in numerical computations of the eigenvalues.

Let $\fats_j^T=(s_{jDD}, s_{jDT}, s_{jE}, s_{d}, s_{de})$ be the
eigenvector of $H_1$ corresponding to $\lambda_j$ and let $s_{jA}=\exp(\mu_{jA}+i\nu_{jA}),\, A=DD, DT, E, d, de$.
The oscillatory eigenvalue $\lambda_1=i\theta_1$ has the eigenvector $\fats_1$.
The eigenvector of $\lambda_2=i\theta_2=\lambda_1^{\ast}=-i\theta_1$ is $\fats_2=\fats_1^\ast$.
Then with initial data $\dfr_0(0)=\fats_1+\fats_1^\ast$, the oscillatory perturbation
is derived from the solution of \eqref{eq:Huang_dfr0}
\begin{equation}
\begin{array}{rll}
   \dr_{0A}(x, t)&=(s_{1A}\exp(i \theta_1 t)+s_{1A}^\ast\exp(-i\theta_1 t))\cos(\pi x/L)\\
                        &=2|s_{1A}|\cos(\nu_{1A}+ \theta_1 t)\cos(\pi x/L),\\
             A&=DD, DT, E, d, de.
\end{array}
\label{eq:Huang_sol}
\end{equation}
The oscillations in time have the period $T=2\pi/\theta_1$ with different phase angles $\nu_{1A}$ for
the species. The period will change when $L$ increases due to cell growth. In the analysis here, we let $L$ be constant.

The spatial mode in \eqref{eq:drhocos} and \eqref{eq:Huang_sol} has two peaks in space, one at $x=0$ and one at $x=L$ with
alternating sign and oscillates in time.
The time average of the square of the species concentration in \eqref{eq:Huang_sol} is
\[
   \frac{1}{T}\int_0^T (\dr_{0A}(x, t))^2\, dt=|s_{1A}|^2\cos^2(\pi x/L)
\]
with a dip at the center of the cell in the MinD concentration as
observed in Figures~\ref{fig:URDMETimeSeries} and \ref{fig:mincdepde}.
Similar analyses for related model equations can be found in \citep{Kruse:2002, MeaKru:2005}.




\begin{figure}[htbp]
\centering
\subfigure{\includegraphics[width=0.8\linewidth]{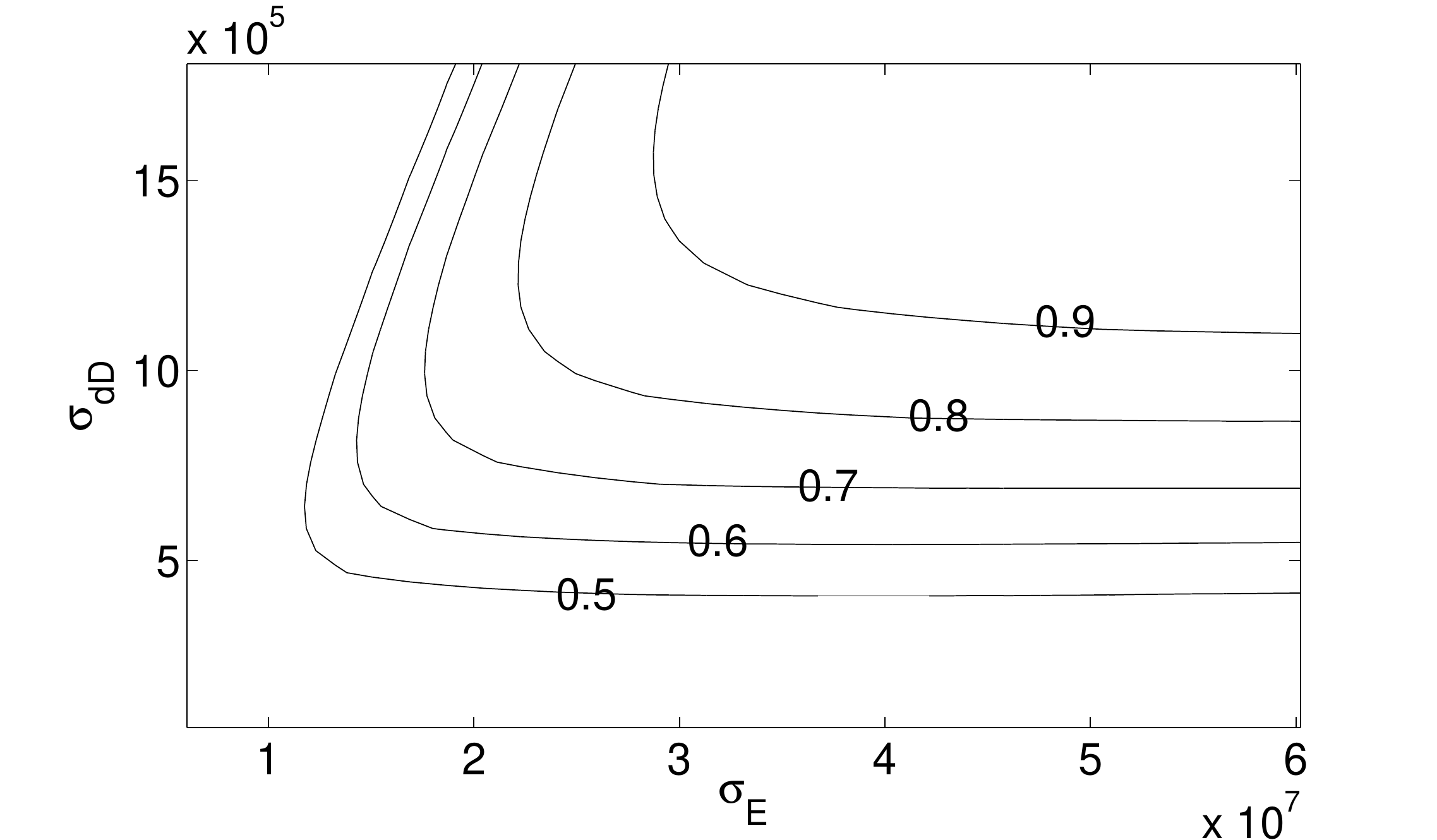}}
\subfigure{\includegraphics[width=0.8\linewidth]{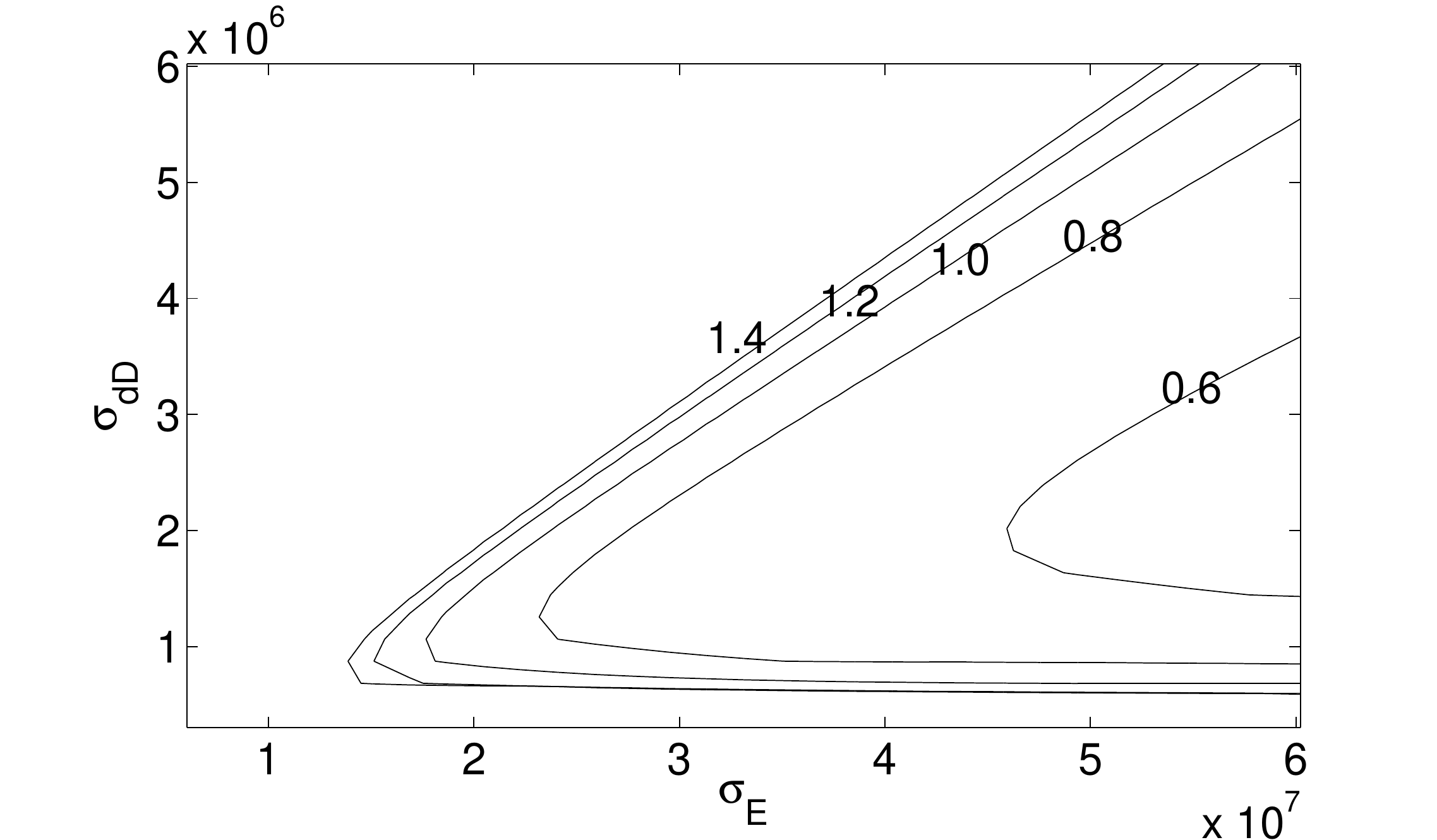}}
\caption{The perturbations of the steady state solution are oscillatory on the lines in the $\sigma_{E}-\sigma_{dD}$ plane. 
The curves represent different values of $\sigma_{de}$ (top) and $\sigma_{DT}$ (bottom).}
\label{fig:perturb}
\end{figure}
The sensitivity in the oscillatory eigenvalue of $H_1$ in \eqref{eq:drhosol}
to changes in the reaction parameters is evaluated in Figure~\ref{fig:perturb}.
The isolines for $\max\Re\lambda_j=0$
are drawn in the $\sigma_{dD}-\sigma_{E}$ plane for different $\sigma_{de}$ and $\sigma_{DT}$.
The eigenvalues are insensitive to $\sigma_D$. This is confirmed in \citep{TouJerRut:2006}.
In the stable regions in the lower and left parts of the figures with $\Re\lambda_j(H_1)<0$, $\dfr_0(x, t)$ will decay in \eqref{eq:drhosol} and
$\fatrho(x, t)$ will approach the steady state $\fatrho_\infty$. The perturbation is mildly unstable in the upper
right part of the figures and will grow there until nonlinear effects, non-negativity and the bounds on
the total number of MinD and MinE molecules \eqref{eq:Huang_cons} will limit the amplitude.
For these parameters and the two oscillatory eigenvalues, $\Re\lambda_j$ is small compared to $\Im\lambda_j$.

\begin{figure}[htbp]
\centering
\includegraphics[width=\linewidth]{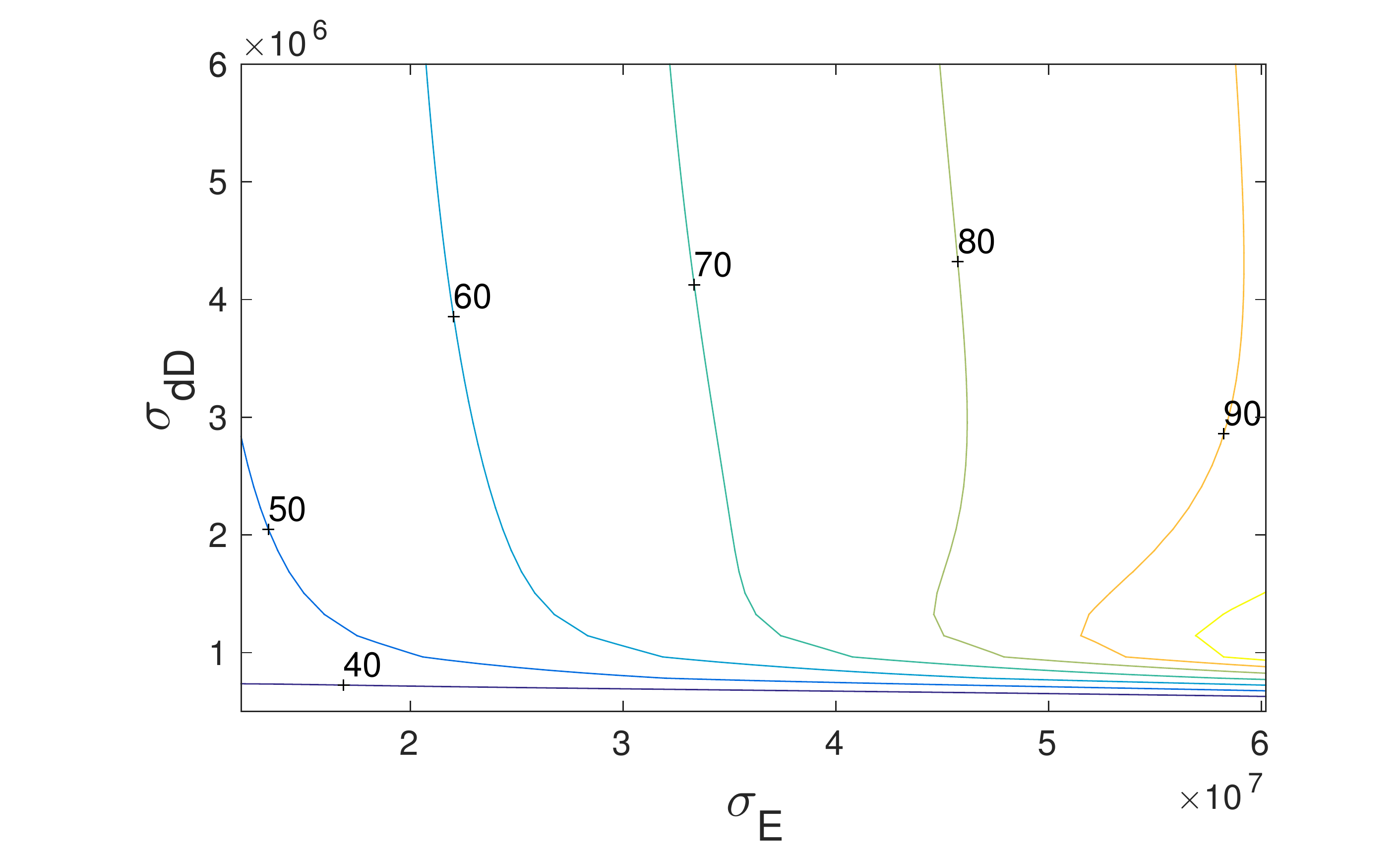}
\caption{The isolines of the period of the oscillations for the perturbations of the steady state
solution in the $\sigma_{E}-\sigma_{dD}$ plane. The curves represent different values of $T$.}
\label{fig:perturb2}
\end{figure}
The dependence of the period on the $\fatsigma$ parameters is displayed in Figure~\ref{fig:perturb2}. The isolines of $T$
are computed as $2\pi/\Im\lambda_j$ for the oscillatory eigenvalues. The period varies quickly when $\sigma_{dD}$ is changed
around the base values of $\fatsigma$ but is insensitive to perturbations in $\sigma_E$ there.


\subsection{Perturbation analysis}\label{sec:perturb}

Equations for the space and time dependent perturbation $\dfr$ of the steady state in \eqref{eq:Huang_pert} and \eqref{eq:drhoeps} when $\veps>0$
will be derived from \eqref{eq:Huang_eps} using separation of variables.

The solution is first expanded in a cosine series in space
\begin{equation}
      \dfr(x, t)=\displaystyle{\sum_{\omega=1}^\infty \hdfr_{\omega}(t)\cos(\omega \pi x/L),\; x\in[0, L],\; t\ge 0}.
\label{eq:dfrcos}
\end{equation}
Then the boundary conditions in \eqref{eq:Huang_1D} are satisfied. The initial condition is taken to be
\begin{equation}
      \dfr(x, 0)=\displaystyle{\hdfr_{1}(0)\cos(\pi x/L)}.
\label{eq:dfrcosini}
\end{equation}
Let $H_{\omega}=J-(\gamma\omega^2\pi^2/L^2) D$ and insert $\dfr$ in \eqref{eq:dfrcos} into \eqref{eq:Huang_eps} to obtain
\begin{equation}
\begin{array}{rll}
      \displaystyle{\pt\sum_{\omega=1}^\infty \hdfr_{\omega}(t)\cos(\omega \pi x/L)}=\displaystyle{\sum_{\omega=1}^\infty H_{\omega}\hdfr_{\omega}\cos(\omega \pi x/L)}\\
              \quad\displaystyle{+\veps(\kappa_t+\kappa_x)\sum_{\omega=1}^\infty G\hdfr_{\omega}\cos(\omega \pi x/L)+\ordo(\veps^2)}.
\end{array}
\label{eq:dfrcoseq}
\end{equation}

Introduce a change of variables $\hdfr_{\omega}=S\hdfu_{\omega}$ where $S=(\fats_{1}, \fats_{2}, \fats_{3}, \fats_{4}, \fats_{5})$
is the eigenvector matrix of $H_1$. The corresponding transformations of $H_{\omega}$ and $G$ are $\hat{H}_\omega=S^{-1}H_{\omega}S$ and $\hG=S^{-1}GS$.
When $\omega=1$, $H_1$ is a diagonal matrix with the eigenvalues $\lambda_j=\lambda_j(H_1)$ on the diagonal. The eigenvalues of $H_{\omega}, \;\omega\ge 2,$
are $\lambda_{\omega j}$. For the linearized system of equations \eqref{eq:Huang_eps}, we {\it assume}
\begin{equation}
\begin{array}{lll}
    D_{ii}\ge 0,\; i=1,\ldots,5,\; D_{jj}>0\; {\rm for}\;{\rm at}\;{\rm least}\;{\rm one}\; j,\\
    \lambda_1=i\theta_1,\;\lambda_2=i\theta_2=-i\theta_1,\; \Re\lambda_j<0,\; j=3,4,5,\\
    \Re\lambda_{\omega j}<0,\; j=1,\ldots,5,\;\omega\ge 2.
\end{array}
\label{eq:lamassump}
\end{equation}
The third assumption concerning $\lambda_{\omega j}$ is not necessary if the diffusion is the same in all components 
with $D=I$. Then $\lambda_{\omega j}=\lambda_j-\gamma(\omega^2-1)\pi^2/L^2$ and 
\begin{equation}
    \Re\lambda_{\omega j}=\Re\lambda_j-\gamma(\omega^2-1)\pi^2/L^2<\Re\lambda_j\le 0,\; j=1,\ldots,5,\;\omega\ge 2.
\label{eq:Bendix}
\end{equation}
The following analysis is also simplified considerably if $D=I$.

The equations satisfied by the coefficients $\hdu_{\omega j}(t)$ are for $j=1,\ldots,5,$
\begin{equation}
\begin{array}{rll}
      \displaystyle{\pt\sum_{\omega=1}^\infty \hdu_{\omega j}\cos(\omega \pi x/L)}=\displaystyle{\sum_{\omega=1}^\infty\sum_{k=1}^5 \hH_{\omega jk}\hdu_{\omega k}\cos(\omega \pi x/L)}\\
              \quad\displaystyle{+\veps(\kappa_t+\kappa_x)\sum_{\omega=1}^\infty\sum_{k=1}^5 \hG_{jk} \hdu_{\omega k}\cos(\omega \pi x/L)+\ordo(\veps^2)}.
\end{array}
\label{eq:ducoseq}
\end{equation}
A Lindstedt-Poincar{\'e} transformation of time
\begin{equation}
   t=s(1+\veps\psi_{j1}(s)+\veps^2\psi_{j2}(s)+\ldots)
\label{eq:ttrans}
\end{equation}
is introduced for the $j$:th equation to avoid secular solutions later with terms in $\hdu_{\omega j}$ growing linearly
in time, see e.g. \citep{Nayfeh}. Then the time derivative is transformed to
\begin{equation}
   \ps\hdu_{\omega j}=\pt\hdu_{\omega j}\frac{dt}{ds}=(1+\veps s\psi_{j1}'+\veps\psi_{j1}+\ordo(\veps^2))\pt\hdu_{\omega j}.
\label{eq:ttransderiv}
\end{equation}
Consequently, the equation in $s$ is
\begin{equation}
\begin{array}{rll}
      \displaystyle{\ps\sum_{\omega=1}^\infty \hdu_{\omega j}\cos(\omega \pi x/L)}=\displaystyle{(1+\veps(s\psi_{j1}'+\psi_{j1}))
                                                                            \sum_{\omega=1}^\infty\sum_{k=1}^5 \hH_{\omega jk}\hdu_{\omega k}\cos(\omega \pi x/L)}\\
              \quad \displaystyle{+\veps(\kappa_t+\kappa_x)\sum_{\omega=1}^\infty \sum_{k=1}^5 \hG_{jk}\hdu_{\omega k}\cos(\omega \pi x/L)+\ordo(\veps^2)}.
\end{array}
\label{eq:ducoseqs}
\end{equation}
Insert the $\veps$-expansion of $\hdu_{\omega j}$
\begin{equation}
   \hdu_{\omega j}=\hdu_{\omega j 0}+\veps\hdu_{\omega j 1}+\ordo(\veps^2)
\label{eq:dueps}
\end{equation}
into \eqref{eq:ducoseqs} and collect terms multiplied by $\veps^k,\; k=0,1,2\ldots$
For $\veps^0$ we arrive at an equation for $\hdu_{\omega j 0}$
\begin{equation}
      \ps\hdu_{\omega j 0}=\sum_{k=1}^5 \hH_{\omega jk}\hdu_{\omega k 0},\; \omega\ge 1.
\label{eq:hdu0eq}
\end{equation}
By \eqref{eq:dfrcosini}, the initial conditions are
\begin{equation}
\begin{array}{rll}
   \hdu_{1j0}(0)&=\hdu_{1j00}= \displaystyle{\sum_{k=1}^5 (S^{-1})_{jk}\hdr_{1 k}(0), \; j=1,\ldots,5},\\
   \hdu_{\omega j 0}(0)&=0,\; \omega\ge 2,\;  j=1,\ldots,5.
\end{array}
\label{eq:hduini}
\end{equation}
By assumption \eqref{eq:lamassump} for $\omega=1$, the solution for large $s$ is
\begin{equation}
\begin{array}{lll}
   \hdu_{110}(s)=\exp(i\theta_1 s)\hdu_{1100},\;\hdu_{120}(s)=\exp(-i\theta_1 s)\hdu_{1200}, \\
   \hdu_{1j0}(x, s)\approx 0,\; j=3,4,5,
\end{array}
\label{eq:hdusol}
\end{equation}
and because of the initial conditions
\begin{equation}
   \hdu_{\omega j 0}(s)=0,\; \omega\ge 2,\;  j=1,\ldots,5.
\label{eq:hdusol2}
\end{equation}

Since $\lambda_j\ne \pm i\theta_1$ when $j=3,4,5,$ there is no secular term for these $j$ and we let $\psi_{ji}=0,\, i\ge 1,$
in \eqref{eq:ttrans} and $s=t$. With the approximations in \eqref{eq:hdusol} and assuming
that $\hdu_{1100}=\hdu_{1200}=1$ to simplify the notation,
the equations for $\hdu_{\omega j1}$ follow from terms proportional to $\veps^1$
\begin{equation}
\begin{array}{lll}
      \displaystyle{\ps\sum_{\omega=1}^{\infty}\hdu_{\omega j 1}\cos(\omega\pi x/L)=\sum_{k=1}^\infty\sum_{k=1}^5\hH_{\omega jk}\hdu_{\omega k1}\cos(\omega\pi x/L)}\\
                                                            \quad   \displaystyle{+(\kappa_t(s)+\kappa_x(x))\cos(\pi x/L)(\hG_{j1}\exp(i\theta_1 s)+\hG_{j2}\exp(-i\theta_1 s))}\\
                \quad  \displaystyle{+(s\psi_{j1}'+\psi_{j1})i\theta_j\exp(i\theta_j s)\cos(\pi x/L),\quad j=1,2,}\\
      \displaystyle{\pt\sum_{\omega=1}^{\infty}\hdu_{\omega j 1}\cos(\omega\pi x/L)=\sum_{k=1}^\infty\sum_{k=1}^5\hH_{\omega jk}\hdu_{\omega k1}\cos(\omega\pi x/L)}\\
                                          \quad  \displaystyle{+(\kappa_t(t)+\kappa_x(x))\cos(\pi x/L)(\hG_{j1}\exp(i\theta_1 t)+\hG_{j2}\exp(-i\theta_1 t)),\quad j=3,4,5.}
\end{array}
\label{eq:du1eq1}
\end{equation}

Let $\kappa_{x}$ have the cosine expansion
\begin{equation}
   \kappa_x(x)=\sum_{\omega=2}^\infty \hat{\kappa}_{x\omega}\cos(\omega \pi x/L)
\label{eq:kappacos}
\end{equation}
such that \eqref{eq:enoise} is fulfilled. The factor $\kappa_x\cos(\pi x/L)$ in \eqref{eq:du1eq1} can be written
\begin{equation}
\begin{array}{ll}
   \kappa_x(x)\cos(\pi x/L)=\displaystyle{\sum_{\omega=1}^\infty \tilde{\kappa}_{x\omega}\cos(\omega \pi x/L)},\\
  \quad \tilde{\kappa}_{x\omega}=\frac{1}{2}\hat{\kappa}_{x,\omega+1},\, \omega=1,2,\quad
   \tilde{\kappa}_{x\omega}=\frac{1}{2}(\hat{\kappa}_{x,\omega-1}+\hat{\kappa}_{x,\omega+1}),\; \omega\ge 3.
\end{array}
\label{eq:kappacos2}
\end{equation}
Using the expansion \eqref{eq:kappacos2} in \eqref{eq:du1eq1} we obtain the equations for $\omega=1$
\begin{equation}
\begin{array}{lll}
\ps\hdu_{1j1}=&i\theta_j\hdu_{1j1}+(\kappa_t+\tilde{\kappa}_{x1})(\hG_{j1}\exp(i\theta_1 s)+\hG_{j2}\exp(-i\theta_1 s))\\
                 &+(s\psi_{j1}'+\psi_{j1})i\theta_j\exp(i\theta_j s),\; j=1,2,\\
\pt\hdu_{1j1}=&\lambda_j\hdu_{1j1}+(\kappa_t+\tilde{\kappa}_{x1})(\hG_{j1}\exp(i\theta_1 t)+\hG_{j2}\exp(-i\theta_1 t)),\; j=3,4,5.
\end{array}
\label{eq:du1omeq}
\end{equation}


Choose $\psi_{j1}(s)$ for $j=1,2,$ in \eqref{eq:ttrans} such that
\begin{equation}
   i\theta_j(s\psi_{j1}'(s)+\psi_{j1}(s))+(\kappa_t(s)+\tilde{\kappa}_{x1})\hG_{jj}=0
\label{eq:Eulereq}
\end{equation}
in \eqref{eq:du1omeq}. Then the equation for $\omega=1$ and $j=1$ is
\begin{equation}
   \ps\hdu_{111}=i\theta_1\hdu_{111}+(\kappa_t+\tilde{\kappa}_{x1})\hG_{12}\exp(-i\theta_1 s)
\label{eq:du1om1eq}
\end{equation}
with the initial condition $\hdu_{111}(0)=0$. The solution to \eqref{eq:du1om1eq} is
\begin{equation}
\begin{array}{rll}
   \hdu_{111}(s)=&\displaystyle{\int_0^s \exp(i\theta_1(s-v))(\kappa_t+\tilde{\kappa}_{x1})\hG_{12}\exp(-i\theta_1 v)\, \ud v}\\
            =&\displaystyle{\hG_{12}\left(\frac{\tilde{\kappa}_{x1}}{\theta_1}\sin(\theta_1 s)+
              \exp(i\theta_1 s)\int_0^s \exp(-2i\theta_1 v)\kappa_t(v)\, \ud v\right).}
\end{array}
\label{eq:du1om1sol1}
\end{equation}
The solution for $j=2$ is obtained by replacing $\theta_1$ by $-\theta_1$ and switching the indices 1 and 2
in \eqref{eq:du1om1sol1}
\begin{equation}
\begin{array}{lll}
   \hdu_{121}(s)=\displaystyle{\hG_{21}\left(-\frac{\tilde{\kappa}_{x1}}{\theta_1}\sin(\theta_1 s)+
              \exp(-i\theta_1 s)\int_0^s \exp(2i\theta_1 v)\kappa_t(v)\, \ud v\right).}
\end{array}
\label{eq:du1om1sol2}
\end{equation}

The solution of \eqref{eq:du1omeq} for $\omega=1$ and $j\ge 3$ is
\begin{equation}
\begin{array}{rll}
   \hdu_{1j1}(t)=&\displaystyle{\int_0^t \exp(\lambda_j(t-v))(\kappa_t+\tilde{\kappa}_{x1})
                (\hG_{j1}\exp(i\theta_1 v)+\hG_{j2}\exp(-i\theta_1 v))\, \ud v}\\
            =&\displaystyle{\frac{\tilde{\kappa}_{x1}\hG_{j1}}{\lambda_j-i\theta_1}
                            (\exp(\lambda_j t)-\exp(i\theta_1 t))}
             \displaystyle{+\frac{\tilde{\kappa}_{x1}\hG_{j2}}{\lambda_j+i\theta_1}
                            (\exp(\lambda_j t)-\exp(-i\theta_1 t))}\\
             &\displaystyle{+\hG_{j1}\exp(\lambda_jt)\int_0^t\kappa_t(v) \exp((-\lambda_j+i\theta_1)v))\, \ud v}\\
             &\displaystyle{+\hG_{j2}\exp(\lambda_jt)\int_0^t\kappa_t(v) \exp(-(\lambda_j-i\theta_1)v))\, \ud v.}
\end{array}
\label{eq:du1om1sol3}
\end{equation}
Since $\Re\lambda_j<0$ by the assumption \eqref{eq:lamassump}, $\exp(\lambda_j t)$ vanishes for large $t$.

The equations for $\omega\ge 2$ and $j=1,\ldots,5,$ are derived from \eqref{eq:du1eq1} and \eqref{eq:kappacos2}
\begin{equation}
\begin{array}{lll}
      \pt\hdu_{\omega j 1}=&\displaystyle{\sum_{k=1}^5 \hH_{\omega jk}\hdu_{\omega k1}
                 +\tilde{\kappa}_{x\omega}(\hG_{j1}\exp(i\theta_1 t)+\hG_{j2}\exp(-i\theta_1 t))}.
\end{array}
\label{eq:duomjeq}
\end{equation}
Transform back in \eqref{eq:duomjeq} from $\hdfu_{\omega}$ to $\hdfr_{\omega}$ in \eqref{eq:dfrcoseq} using $S$.
The eigenvector matrix of $H_{\omega}$ is
$S_{\omega}=(\fats_{\omega 1}, \fats_{\omega 2}, \fats_{\omega 3}, \fats_{\omega 4}, \fats_{\omega 5})$ and the eigenvalues $\lambda_{\omega j}$ satisfy $\Re\lambda_{\omega j}<0$ by \eqref{eq:lamassump}.
Then change the variables such that $\hdfr_{\omega}=S_{\omega}\tdfu_{\omega}$. Let $\tdfu_{\omega 1}$ be the term in $\tdfu_{\omega}$ multiplied by $\veps$.
The equation for the $j$:th component of $\tdfu_{\omega 1}$ is
\begin{equation}
\begin{array}{lll}
      \pt\tdu_{\omega j 1}=&\lambda_{\omega j}\tdu_{\omega j 1}
                 +\tilde{\kappa}_{x\omega}(\tG_{\omega j1}\exp(i\theta_1 t)+\tG_{\omega j2}\exp(-i\theta_1 t)),
\end{array}
\label{eq:duomjeq2}
\end{equation}
where $\tG_{\omega j\ell}=\sum_{k=1}^5 (S_{\omega}^{-1}S)_{jk}\hG_{k\ell}, \;\ell=1,2$.
Solving \eqref{eq:duomjeq2} for $\omega\ge 2$ using the initial conditions \eqref{eq:hduini} we arrive at a solution similar to \eqref{eq:du1om1sol3}
\begin{equation}
\begin{array}{rll}
   \tdu_{\omega j1}(t)=&\displaystyle{\int_0^t \exp(\lambda_{\omega j}(t-v))\tilde{\kappa}_{x\omega}
                (\tG_{\omega j1}\exp(i\theta_1 v)+\tG_{\omega j2}\exp(-i\theta_1 v))\, \ud v}\\
            =&\displaystyle{\frac{\tilde{\kappa}_{x\omega}\tG_{\omega j1}}{\lambda_{\omega j}-i\theta_1}
                  (\exp(\lambda_{\omega j} t)-\exp(i\theta_1 t))}\\
             &\displaystyle{+\frac{\tilde{\kappa}_{x\omega}\tG_{\omega j2}}{\lambda_{\omega j}+i\theta_1}
                   (\exp(\lambda_{\omega j} t)-\exp(-i\theta_1 t))}.
\end{array}
\label{eq:du1om2sol}
\end{equation}
For large $t$, $\exp(\lambda_{\omega j} t)\rightarrow 0$ and $\hdu_{\omega j1},\, \omega=1,2,\ldots$ is simplified to
\begin{equation}
   \tdu_{\omega j1}(t)=-\tilde{\kappa}_{x\omega}\Theta_{\omega j},\;\Theta_{\omega j}=\displaystyle{\frac{\tG_{\omega j1}}{\lambda_{\omega j}-i\theta_1}
                  \exp(i\theta_1 t)+\frac{\tG_{\omega j2}}{\lambda_{\omega j}+i\theta_1}
                   \exp(-i\theta_1 t)}.
\label{eq:du1om2solsim}
\end{equation}

The solution to the Euler equation \eqref{eq:Eulereq} is
\begin{equation}
   \psi_{j1}(s)=-\frac{\hG_{jj}}{i\theta_j s}\left(\tilde{\kappa}_{x1} s+\int_0^s \kappa_t(v)\, \ud v\right).
\label{eq:Eulersol}
\end{equation}
Introduce
\begin{equation}
   \kappa_{xt}(t)=\tilde{\kappa}_{x1} +\frac{1}{t}\int_0^t \kappa_t(v)\, \ud v.
\label{eq:kappaxtdef}
\end{equation}
By \eqref{eq:ttrans} and \eqref{eq:Eulersol} we conclude that
\begin{equation}
   s=\displaystyle{t\left(1+\frac{\veps \hG_{jj}}{i\theta_j}\kappa_{xt}(s)\right)+\ordo(\veps^2)
    =t\left(1+\frac{\veps \hG_{jj}}{i\theta_j}\kappa_{xt}(t)\right)+\ordo(\veps^2),\; j=1,2}.
\label{eq:invttrans}
\end{equation}

Let $\tdfu_{\omega}$ with components $\tdu_{\omega j}$ be defined
as in \eqref{eq:dueps}. Combining \eqref{eq:dueps}, \eqref{eq:hdusol}, \eqref{eq:du1om1sol1}, and \eqref{eq:du1om2solsim} with \eqref{eq:invttrans},
the two lowest order terms in the $\veps$-expansion of $\hdu$ are
\begin{equation}
\begin{array}{lll}
   \hdu_{11}(t)=\displaystyle{\hdu_{110}(t)+\veps\hdu_{111}(t)+\ordo(\veps^2)=\exp\left(i\theta_1 t+\veps \hG_{11}t\kappa_{xt}(t)\right)}\\
           \quad +\displaystyle{\veps\hG_{12}\left(\frac{\tilde{\kappa}_{x1}}{\theta_1}\sin(\theta_1 t)+
              \exp(i\theta_1 t)\int_0^t \exp(-2i\theta_1 v)\kappa_t(v)\, \ud v\right)+\ordo(\veps^2)}\\
   \hdu_{1j}(t)= \displaystyle{\veps\hdu_{1j1}(t)+\ordo(\veps^2)}=
       \displaystyle{-\veps\tilde{\kappa}_{x1}\left(\frac{\hG_{j1}}{\lambda_j-i\theta_1}\exp(i\theta_1 t)+\frac{\hG_{j2}}{\lambda_j+i\theta_1}\exp(-i\theta_1 t)\right)}\\
          \quad    \displaystyle{+\veps\hG_{j1}\exp(\lambda_jt)\int_0^t\kappa_t(v) \exp((-\lambda_j+i\theta_1)v))\, \ud v}\\
         \quad    \displaystyle{+\veps\hG_{j2}\exp(\lambda_jt)\int_0^t\kappa_t(v) \exp(-(\lambda_j-i\theta_1)v)\, \ud v+\ordo(\veps^2),\quad j=3,4,5,}\\
   \tdu_{\omega j}(t)=\displaystyle{\veps\tdu_{\omega j1}(t)+\ordo(\veps^2)}
         \quad =\displaystyle{-\veps \tilde{\kappa}_{x\omega}\Theta_{\omega j}(t)+\ordo(\veps^2),\; \omega\ge 2,}
\end{array}
\label{eq:du1sol}
\end{equation}
when the transient has disappeared for large $t$.
The solution for the other oscillatory mode $\hdu_{12}$ is obtained from $\hdu_{11}$ by replacing $\theta_1$ by $-\theta_1$
and switching the indices $1\rightarrow 2$ and $2\rightarrow 1$ in $\hG_{jk}$ as in \eqref{eq:du1om1sol2}.
If $\kappa_x=0$ then \eqref{eq:du1sol} is simplified and $\tdu_{\omega j}(t)=\ordo(\veps^2)$.

In the original variables, we have from \eqref{eq:du1sol} that
\begin{equation}
\begin{array}{lll}
   \dfr&=\displaystyle{S_1\hdfu_1\cos(\pi x/L)+\sum_{\omega=2}^\infty S_{\omega}\tdfu_{\omega}\cos(\omega\pi x/L)}\\
       &=\displaystyle{(\hdu_{110}(t)\fats_{11}+\hdu_{120}(t)\fats_{12})\cos(\pi x/L)+\ordo(\veps)}.
\end{array}
\label{eq:drhoexpr}
\end{equation}
The main oscillatory mode given by $\fats_{11}$ and $\fats_{12}$ is perturbed by a term of $\ordo(\veps)$ due to the
perturbed coefficients in the model.

The inverse of the eigenvector matrix $S_{\omega}^{-1},\; \omega\ge 1,$ has the properties
\begin{equation}
   S_{\omega}^{-1}=\left(\begin{array}{c} \fatvarsigma_{\omega 1}\\ \fatvarsigma_{\omega 2}\\ \ldots\\\fatvarsigma_{\omega 5}
         \end{array}\right),\;
   S_{\omega}^{-1} S_{\omega} = \left(\begin{array}{cccc} \fatvarsigma_{\omega 1} \fats_{\omega 1}& \fatvarsigma_{\omega 1} \fats_{\omega 1}^\ast& \ldots\\
                                       \fatvarsigma_{\omega 2} \fats_{\omega 1}& \fatvarsigma_{\omega 2} \fats_{\omega 1}^\ast& \ldots\\
                                        \ldots
         \end{array}\right)
             = \left(\begin{array}{cccc} 1& 0& \ldots\\
                                         0& 1& \ldots\\
                                        \ldots
         \end{array}\right).
\label{eq:invprop}
\end{equation}
The rows $\fatvarsigma_{\omega j}$ of $S_{\omega}^{-1}$ are such that $\fatvarsigma_{\omega j} \fats_{\omega k}=0$ when $j\ne k$ and $\fatvarsigma_{\omega j} \fats_{\omega j}=1$.
We find that $\fatvarsigma_{\omega 2}=\fatvarsigma_{\omega 1}^\ast$ by $S_{\omega}^{-1} S_{\omega}$ in \eqref{eq:invprop}. The elements in the upper left corner
of $\tG_\omega$ and $\hG=\tG_1$ are then
\begin{equation}
   \hG_{\omega 11}=\fatvarsigma_{\omega 1} G \fats_{\omega 1}=\fatvarsigma_{\omega 2}^\ast G \fats_{\omega 2}^\ast=\hG_{\omega 22}^\ast,\;
   \hG_{\omega 12}=\fatvarsigma_{\omega 1} G \fats_{\omega 2}=\fatvarsigma_{\omega 2}^\ast G \fats_{\omega 1}^\ast=\hG_{\omega 21}^\ast.
\label{eq:Gexpl}
\end{equation}
Therefore, $\hdu_{120}=\hdu_{110}^\ast$ in \eqref{eq:du1sol} and $\dfr$ in \eqref{eq:drhoexpr} is
\begin{equation}
\begin{array}{lll}
   \dfr&=(\hdu_{110}\fats_1+\hdu_{110}^\ast\fats_1^\ast)\cos(\pi x/L)+\ordo(\veps)\\
       &=2 \Re \{\hdu_{110}(t)\fats_1\}\cos(\pi x/L)+\ordo(\veps).
\end{array}
\label{eq:drhoexpr2}
\end{equation}

The argument in the exponential in the leading term of $\ordo(1)$ in $\hdu_{110}$ in \eqref{eq:du1sol} is denoted by $(\xi+i\eta)t$ with
\begin{equation}
     \xi(t)=\veps\hG_{11R}\kappa_{xt}(t),\quad \eta(t)= \theta_1 +\veps\hG_{11I}\kappa_{xt}(t),\quad
       \hG_{11R}=\Re\hG_{11},\;  \hG_{11I}=\Im\hG_{11}.
\label{eq:xieta}
\end{equation}
Depending on the sign of $\xi(t)$ in \eqref{eq:xieta}, there will be a slow growth or decay of the main oscillatory mode.
The frequency of the oscillations in $\eta(t)$ will be perturbed slightly depending on $\kappa_{xt}(t)$. Thus, by \eqref{eq:drhoexpr2}
and as in \eqref{eq:Huang_sol}
\begin{equation}
\begin{array}{rll}
   \dr_A=&2\Re \{\exp\left(i(\theta_1 +\veps \hG_{11I}\kappa_{xt}(t))t\right)\exp\left(\veps \hG_{11R}t\kappa_{xt}(t)\right)\\
          &\times\cos(\pi x/L)|s_{1A}|\exp(i\nu_{1A})\}+\ordo(\veps)\\
       =&2|s_{1A}|\exp\left(\veps \hG_{11R}t\kappa_{xt}(t)\right)\cos\left((\theta_1 +\veps \hG_{11I}\kappa_{xt}(t))t+\nu_{1A}\right)\\
          &\times\cos(\pi x/L)+\ordo(\veps),\\
      A=&DD, DT, E, de, e.
\end{array}
\label{eq:drhoexpr3}
\end{equation}
The oscillations in all components are modified by $\veps \hG_{11I}\kappa_{xt}(t)$ in \eqref{eq:drhoexpr3} due to
an accumulation of the temporal perturbation $\kappa_t$ in \eqref{eq:kappaxtdef}
and the constant spatial perturbation $\tilde{\kappa}_{x1}=\frac{1}{2}\hat{\kappa}_{x2}$ in
\eqref{eq:kappacos2} and \eqref{eq:kappacos} of the reaction coefficients. The temporal perturbation vanishes for large $t$
by the assumption in \eqref{eq:enoise}.

Partition the interval $[0, t]$ into subintervals $[t_{j-1}, t_j], j=1,\ldots, J,$ with $\Delta t_j=t_j-t_{j-1}$
and use \eqref{eq:kappaxtdef}.
The dominant part of $\hdu_{11}$ in \eqref{eq:du1sol} will evolve between $t_{j-1}$ and $t_j$ as
\begin{equation}
   \hdu_{110}(t_j)=\displaystyle{\exp\left(i\theta_1 \Delta t_j+\veps \hG_{11}(\tilde{\kappa}_{x1}\Delta t_j
                +\int_{t_{j-1}}^{t_j} \kappa_t(v)\, \ud v)\right)\hdu_{110}(t_{j-1})+\ordo(\veps)}.
\label{eq:du1sol2}
\end{equation}
Introduce the average
\begin{equation}
     \kappa_{tj}=\frac{1}{\Delta t_j}\int_{t_{j-1}}^{t_j} \kappa_t(v)\, \ud v
\label{eq:kappatj}
\end{equation}
in \eqref{eq:du1sol2}. Then
\begin{equation}
\begin{array}{rll}
   \hdu_{110}(t_j)&=\displaystyle{\exp\left(i\theta_1 \Delta t_j+\veps \hG_{11}\Delta t_j(\tilde{\kappa}_{x1}
                +\kappa_{tj})\right)\hdu_{110}(t_{j-1})+\ordo(\veps)}\\
              &=\displaystyle{\exp\left(i\theta_1 t_j+\veps \hG_{11} (t_j\tilde{\kappa}_{x1}
                +\sum_{k=1}^J\Delta t_j \kappa_{tj})\right)\hdu_{110}(0)+\ordo(\veps)},
\end{array}
\label{eq:du1sol3}
\end{equation}
and the frequency $\eta(t_j)$ in \eqref{eq:xieta} at $t_j$ is
\begin{equation}
     \eta(t_j)= \theta_1 +\veps\hG_{11I}(\tilde{\kappa}_{x1}
                +\frac{1}{t_j}\sum_{k=1}^J \Delta t_k\kappa_{tk})
              = \theta_1 +\veps\hG_{11I}(\tilde{\kappa}_{x1}
                +\frac{1}{t_j}\int_0^{t_j}\kappa_{t}(v)\, \ud v).
\label{eq:etatj}
\end{equation}
The contribution to the oscillation in $\hdu_{110}(t_j)$ locally in $[t_{j-1}, t_j]$ is $i\eta_j\Delta t_j$ with
\begin{equation}
     \eta_{j}=\theta_1 +\veps\hG_{11I}(\tilde{\kappa}_{x1}+\kappa_{tj}).
\label{eq:etaj}
\end{equation}

The effects of the perturbations evaluated numerically in Figures~\ref{fig:perturb} and \ref{fig:perturb2} are compared to  
\eqref{eq:drhoexpr3}. The $\fatsigma$ values in 
the unperturbed $\fatf$ in \eqref{eq:Huang_1D} and its Jacobian $F$ are as in \eqref{eq:Huang_param} with the frequency $\theta_1\approx 2\pi/40$. 
Perturbations are introduced 
in each one of the $\fatsigma$ parameters keeping the other ones constant. This defines $G$ in \eqref{eq:Huang_eps}. Then $\hG_{11}$ is computed with the results in Table~\ref{tab:Gperturb}.
\begin{table}[htbp]
\centering
\begin{tabular}{|l|c|c|c|c|c|}
\hline
$\sigma_A$ &$\sigma_{DT}$ &$\sigma_{de}$&$\sigma_{D}$  &$\sigma_{dD}$  &$\sigma_{E}$ \\
\hline
$\hG_{11R}$ & 0.0005 & -0.0654 & -0.0022 &  0.1327 & 0.0074 \\
$\hG_{11I}$ & 0.0343 & 0.2670  & 0.0021  &  0.1999 & 0.0611 \\
\hline
\end{tabular}
\caption{Perturbations in $\hG_{11}$ due to perturbations in the $\fatsigma$ parameters of the model.}
\label{tab:Gperturb}
\end{table}
If $\sigma_{dD}$ is increased then the frequency of the oscillations increases since $\hG_{11I}>0$ in the table and the amplitude
of the linearization increases since $\hG_{11R}>0$ in agreement with Figures~\ref{fig:perturb} and \ref{fig:perturb2}. The sensitivity to
perturbations in $\sigma_{dD}$ is also large in the stochastic simulations in Figure~\ref{fig:URDMETimeSeries}. When $\sigma_{de}$
increases in the left panel of Figure~\ref{fig:perturb} the oscillations are damped and $\hG_{11R}<0$ in the table. The changes 
in stability are small in the right panel of Figure~\ref{fig:perturb} when $\sigma_{DT}$ is varied because $\hG_{11R}\approx 0$.
Neither the frequency nor the amplitude is sensitive to changes in $\sigma_{D}$ in Section~\ref{sec:stabanal} and in Table~\ref{tab:Gperturb}. 
A perturbation in $\sigma_{E}$ has little influence on the frequency in Figure~\ref{fig:perturb2} and $\hG_{11I}$ is small.

Let us consider a special problem where all coefficients are perturbed in the same way
by $\kappa(x, t)$ and $G$ is such that it is diagonalized by $S$. Then $\hG_{jk}=0$ when $j\ne k$ and
the expressions in \eqref{eq:du1sol} are somewhat simplified with $\hdu_{1j}=\hdu_{1j1}=0$ in \eqref{eq:du1om1sol3} for $j\ge 3$.
With this assumption in \eqref{eq:du1sol}
\begin{equation}
   \hdu_{11}(t)=\displaystyle{\exp\left(i\theta_1 t+\veps \hG_{11}t\kappa_{xt}(t))\right)
           +\ordo(\veps^2)}.
\label{eq:du1solspec}
\end{equation}
The effect on $\dfr$ of the higher order spatial modes with $\omega\ge 2$ in $\hdu_{\omega j}$ from \eqref{eq:du1sol} is 
\begin{equation}
    \veps\fath(x, t)=\displaystyle{\sum_{\omega=2}^\infty\sum_{j=1}^5\fats_{\omega j}\tdu_{\omega j}\cos(\omega\pi x/L)}
           =\displaystyle{-\veps\sum_{\omega=2}^\infty\tilde{\kappa}_{x\omega}\cos(\omega\pi x/L)\sum_{j=1}^5\fats_{\omega j}\Theta_{\omega j}(t)}.
\label{eq:du1solspec2}
\end{equation}
The perturbation $\veps\fath$ depends on the spatial perturbations in $\tilde{\kappa}_{x \omega}$ and $\Theta_{\omega j}(t)$ in \eqref{eq:du1om2solsim} which is time dependent
but is independent of $\kappa_x$ and $\kappa_t$.
The expression for $\dfr$ is obtained as in \eqref{eq:drhoexpr2}
\begin{equation}
   \dfr=\du_{1}\fats_1+\du_{1}^\ast\fats_1^\ast+\veps\fath(x, t)+\ordo(\veps^2)=2\Re \{\du_{1}\fats_1\}+\veps\fath(x, t)+\ordo(\veps^2).
\label{eq:drhoexpr4}
\end{equation}
The components of $\dfr$ are
\begin{equation}
\begin{array}{rll}
   \dr_A(x, t)&=&2\Re\{\exp(i\theta_1t+\veps (\hG_{11R}+i\hG_{11I}) t\kappa_{xt}(t))|s_{1A}|\exp(i\nu_{1A})\}\\
              &&   \times \cos(\pi x/L) +\veps h_A(x, t)+\ordo(\veps^2)\\
        &=&2|s_{1A}|\exp(\veps \hG_{11R}t\kappa_{xt}(t))\cos\left((1 +\veps\hG_{11I}\theta_1^{-1}\kappa_{xt}(t))\theta_1t+\nu_{1A}\right)\\
       &&\times \cos(\pi x/L)+\veps h_A(x, t)+\ordo(\veps^2),\\
      A&=&DD, DT, E, d, de.
\end{array}
\label{eq:drhoexpr5}
\end{equation}

The amplitude of the oscillatory mode with frequency $\theta_1$
is perturbed in \eqref{eq:drhoexpr5} by the temporal perturbations $\kappa_t$ if $\Re\hG_{11}\ne 0$ and by the spatial perturbations in $\fath(x)$ and ${\kappa}_{x}$.

\section{Extrinsic noise}\label{sec:extnoise}
We now describe the model for fluctuations in the rate constants.
Briefly, the usual rate constant, $\sigma$, is replaced by $\sigma \kappa(x,t)$, where $\kappa$ is an independent stochastic process.
The fluctuations are such that on average the rate constant is not changed, i.e. $\sigma \kappa(x,t)$ averages to $\sigma$.

\subsection{The Ornstein-Uhlenbeck process}\label{sec:OU}

An Ornstein-Uhlenbeck (OU) process \citep{Gil92,Gil96,Oks98,Kam01} is a scalar, continuous-time, continuous-state Markov process $X(t)$,
that satisfies the stochastic differential equation (SDE)
\begin{equation}
  \ud X = - \frac{1}{\tau} X(t) \ud t + \sqrt{c}\, \ud W .
 \label{eq:OU}
 \end{equation}
Here $W(t)$ is a Wiener process. It is an almost surely continuous function, with $W(0)=0$, increments
$  W(t) - W(s) = \Delta W$ that have a normal distribution with mean $0$ and variance $|t-s|$, and that are
independent on nonoverlapping time intervals.
The positive constant $\tau>0$ is the relaxation time, which is a measure of the average time it
takes the OU process to revert back to the long term mean of $0$ after a fluctuation away from $0$.
For example, the autocorrelation of the OU process is $e^{-t/\tau}$. The parameter $\tau$ is the time scale for how $X(t)$ is correlated in time.
The explicit solution of this SDE (\ref{eq:OU}), and thus a sample path representation of an OU process,
is
\[
  X(t) = e^{-{t}/{\tau}} \big( X(0) + \int_0^{t} \sqrt{c} \, e^{{s}/{\tau}}\ud W(s) \big).
\]
This is one way to see that, given a sure initial condition, $X(0) = x_0$, $X(t)$ has a normal distribution
with mean $x_0 e^{-{t}/{\tau}}$, and variance $\frac{c \tau}{2} (1-e^{-{2t}/{\tau}})$.

More generally, the associated PDF $p(x,t)$ evolves according to the Fokker-Planck PDE
\[
 \pt p(x,t) = \frac{1}{\tau}   \px \big( x p(x,t) \big)
  + \frac{c}{2}  \px^2 p(x,t)  ,
\]
given an initial distribution $p(x,0)$ for $X(0)$.
The diffusion constant $c>0$ controls the spread of the distribution.
For example, as $t \rightarrow \infty$ the process tends to the \textit{stationary distribution}
$p_\infty(x)$, which is normally distributed $\mathcal{N}(0, \frac{c\tau}{2})$ with mean $0$ and variance $\frac{c \tau}{2}$.
If we choose the initial distribution to be the same as the stationary distribution, then the distribution
of $X(t)$ is always the same: $p(x,0) \equiv p_{\infty}(x) = p(x,t) = \mathcal{N}(0, \frac{c\tau}{2})$.
One nice property of the OU process is that it is \textit{ergodic}, i.e. for a sure initial condition, $X(0) = x_0$, and a suitably smooth function $f$
we have
\begin{equation}
 \lim_{t \rightarrow \infty} \frac{1}{t} \int_0^T f \big( X(s) \big) \ud s =  \int_{-\infty}^{\infty} f \big( X \big) \ud \mathbb{P}(X)\equiv \mathbb{E}\left( f(X) \right)=\langle\; f(X)\; \rangle ,
 \label{eq:ergodic}
 \end{equation}
where the second integral is with respect to the probability measure $\mathbb{P}$ that is equal to the stationary distribution $p_\infty(x)$
and $\langle\; \cdot\; \rangle$ denotes the mean value.

To mathematically model the effects of extrinsic noise, we replace the rate constant $\sigma$ by $Y(t)\sigma$.
The extrinsic noise process $Y(t)$ is modeled, up to a normalization constant, by the exponential of the OU process $Y(t) \propto e^{X(t)}$.
One reason for this choice of an exponential is that it ensures that the extrinsic process and the rate $Y(t) \sigma$ are
always positive.
We model $X(t)$ in (\ref{eq:OU}) as always at the stationary distribution, which is $\mathcal{N}(0, \frac{c\tau}{2})$.
Then the \textit{autocorrelation} of $e^{X(t)}$ is
\begin{equation}
 A_{\exp(X)}(s, t)=\frac{\mathbb{E}(e^{X(s)} e^{X(t)}) -\mathbb{E}(e^{X})^2}{{\rm Var}(e^{X})}
          = \frac{\exp({\frac{c \tau}{2} e^{-\frac{|s-t|}{\tau}}}) - 1 }{e^{\frac{c \tau}{2}} - 1}  .
 \label{eq:autocorrelation:expou}
\end{equation}

Notice that $Y$ has a \textit{lognormal distribution} because $X$ has a normal distribution.
When a normal distribution has mean $\mu$ and variance $\sigma^2$, then the mean of the corresponding
lognormal distribution is $e^{\mu + \sigma^2/2}$, and the variance is $(e^{\sigma^2} - 1) e^{2 \mu + \sigma^2} $.
Since the stationary distribution for $X$ is $\mathcal{N}(0, \frac{c\tau}{2})$,  $\mathbb{E}(e^{X}) = e^{c \tau / 4}$ and
the variance ${\rm Var}(e^{X})$ is $e^{\frac{c \tau}{2}}(e^{\frac{c \tau}{2}}-1)$.
We normalize so that the extrinsic noise process has mean $1$ by letting 
\[
Y(t) = e^{X(t)} / \mathbb{E}(e^{X}).
\]
Now with $f(x) = e^{x}$ in (\ref{eq:ergodic}), the ergodic property tells us that the long time
average of the extrinsic noise process $Y(t)$ is $1$.
Thus, by introducing extrinsic noise in this way, on average, we do not change the original value of the rate constant
$\mathbb{E}(Y(t)\sigma) = \sigma \mathbb{E}(Y) = \sigma$.

\subsection{Spatially correlated noise}\label{sec:spatnoise}

The random perturbations at discrete points in space $\veps\kappa_x(x_i),\, i=1,\ldots,N,$ are sampled from a multivariate normal distribution.
The mean value of the perturbations is 0.
The elements of the symmetric, positive definite covariance matrix $C$ are $C_{ij}=\mathbb{E}(\kappa_x(x_i)\kappa_x(x_j)),\, i,j=1,\ldots,N$.
The perturbations are generated by multiplying a vector with independent, normally distributed $\mathcal{N}(0, 1)$ components by the Cholesky factorization of $C$.

%

\section{Random perturbations in time and space}\label{sec:randpert}

The perturbations $\kappa_t(t)$ and $\kappa_x(x)$ in time and space of the parameters in the model in \eqref{eq:Huang_model}
are assumed to be random due to extrinsic noise.
The equation \eqref{eq:Huang_1D} is then a random differential equation with parameters depending on the realization of the process. The analysis of the effect of the perturbations
in Section~\ref{sec:perturb} is the same for deterministic and stochastic perturbations but additional conclusions can be drawn from the distribution
of the stochastic perturbations. As in Section~\ref{sec:extnoise}, the temporal extrinsic noise $\kappa_t(t)$ is here assumed to be generated by an
OU process as in \citep{ShaOll08},
and the spatial extrinsic noise $\kappa_x(x)$ by a multivariate normal distribution.

\subsection{Random perturbations in time}\label{sec:randtime}

Assume that $\kappa_t(t)$ in \eqref{eq:kappadef} and \eqref{eq:kappaxtdef} is a random variable generated by
an OU process $Y(t)$ as in Section \ref{sec:OU} and let $\kappa_t(t)=Y(t)-1$. Then according to \eqref{eq:ergodic}
\begin{equation}
     \lim_{T\rightarrow\infty}\frac{1}{T}\int_{0}^T \kappa_t(v)\, \ud v=\lim_{T\rightarrow\infty}\frac{1}{T}\int_{0}^T Y(v)-1\, \ud v=0.
\label{eq:kappainf}
\end{equation}
The frequency in \eqref{eq:xieta} is
\begin{equation}
     \eta(t)=\theta_1 + \delta_x +\delta_t(t),\quad \delta_x=\veps\hG_{11I}\tilde{\kappa}_{x1}=\frac{1}{2}\veps\hG_{11I}\hat{\kappa}_{x2},\;
        \delta_t(t)=\veps\hG_{11I}\frac{1}{t}\int_{0}^t \kappa_t(v)\, \ud v.
\label{eq:deltadef}
\end{equation}
The deviations from $\theta_1$ caused by the spatial and temporal perturbations are $\delta_x$ and $\delta_t(t)$, respectively.
It follows from \eqref{eq:kappainf} and \eqref{eq:deltadef} that $\lim_{t\rightarrow\infty}\delta_t(t)=0$ as required in \eqref{eq:enoise}.

Define as in \eqref{eq:du1sol2} and \eqref{eq:kappatj}
\begin{equation}
     \eta_j=\theta_1 + \delta_x +\delta_{t_j},\quad
        \delta_{t_j}=\veps\hG_{11I}\kappa_{tj}=\veps\hG_{11I}\frac{1}{\Delta t_j}\int_{t_{j-1}}^{t_j} \kappa_t(v)\, \ud v,
\label{eq:deltadefj}
\end{equation}
in a time interval $[t_{j-1}, t_j]$. In a short interval, we have the approximate instantaneous value $\eta_j$ of the frequency with
the contribution
\begin{equation}
      \delta_{t_j}\approx \veps\hG_{11I}\kappa_t(t_j)
\label{eq:deltadefj2}
\end{equation}
from the variation of $\fatsigma$ in time. The autocorrelation for $\delta_{t_j}$ with a small $\Delta t_j$ is by \eqref{eq:autocorrelation:expou}
\begin{equation}
\begin{array}{ll}
   {A}_t(t_j, t_k)&=\displaystyle{\langle\delta_{t_j}\delta_{t_k}\rangle=
   \frac{\veps^2\hG_{11I}^2}{\Delta t_j \Delta t_k}\int_{t_{j-1}}^{t_j} \int_{t_{k-1}}^{t_k} \langle\kappa_t(u) \kappa_t(v)\rangle\, \ud u\, \ud v}\\
   &\displaystyle{\approx \veps^2\hG_{11I}^2 \langle\kappa_t(t_j) \kappa_t(t_k)\rangle= \veps^2\hG_{11I}^2
    \left( \exp \left( \frac{c \tau}{2} e^{- \frac{\lvert t_j-t_k \rvert}{\tau}} \right) - 1\right)}.
\end{array}
\label{eq:autocorrtj}
\end{equation}
When $\Delta t_{jk}=\lvert t_j-t_k \rvert$ is small compared to $\tau$ in \eqref{eq:autocorrtj} then
\begin{equation}\label{eq:autocorrtappr1}
   \exp \left( \frac{c \tau}{2} e^{- \frac{\Delta t_{jk}}{\tau}} \right) - 1\approx \exp\left(\frac{c \tau}{2}\right)\left(1-\frac{c\Delta t_{jk}}{2}\right)-1,
\end{equation}
and when it is large
\begin{equation}\label{eq:autocorrtappr2}
   \exp \left( \frac{c \tau}{2} e^{- \frac{\Delta t_{jk}}{\tau}} \right) - 1\approx \frac{c \tau}{2}\exp\left(-\frac{\Delta t_{jk}}{\tau}\right).
\end{equation}

\subsection{Random perturbations in space}\label{sec:randspace}

Assume that the spatial perturbation $\kappa_x$ has the cosine expansion with nonzero even coefficients in \eqref{eq:kappacos}
\begin{equation}
   \kappa_x(x)=\sum_{\omega=2,4,\ldots}^\infty \hat{\kappa}_{x\omega}\cos(\omega \pi x/L)
              =\sum_{\mu=1}^\infty \hat{\kappa}_{x,2\mu}\cos(2\mu \pi x/L),\quad x\in[0, L].
\label{eq:kappasumeven}
\end{equation}
Then $\kappa_x(x)$ is $L$-periodic, $\kappa_x(x+L)=\kappa_x(x)$.
The correlation function in space $a_x(\xi)$ is defined by
\begin{equation}
\begin{array}{ll}
   a_x(\xi)&=\displaystyle{\frac{1}{L}\int_0^L \kappa_x(x+\xi)\kappa_x(x)\, \ud x}\\
        &=\displaystyle{\frac{1}{L}\sum_{\mu=1}^\infty\sum_{\nu=1}^\infty\hat{\kappa}_{x,2\mu}^2
        \int_0^L \cos(2\mu\pi(x+\xi)/{L})\cos(2\nu \pi{x}/{L})\, \ud x}
        ,\quad \xi\in[0, L].
\end{array}
\label{eq:corrx}
\end{equation}
This function is also $L$-periodic. Since
\[
    \int_0^L\cos(2\mu\pi x/{L})\cos(2\nu\pi x/{L})\, \ud x=\left\{\begin{array}{rl}
                                                                       0,& \mu\ne \nu,\\
                                                                       L/2,&\mu=\nu,
                                                                     \end{array}\right.
\]
$a_x(\xi)$ can be written
\begin{equation}
   a_x(\xi)=\frac{1}{2}\sum_{\mu=1}^\infty\hat{\kappa}_{x,2\mu}^2 \cos(2\mu\pi\xi/{L})
        ,\quad \xi\in[0, L].
\label{eq:corrx2}
\end{equation}
The coefficients in the cosine expansion of $a_x(\xi)$ are $\hat{\kappa}_{x,2\mu}^2/2$.

When $\kappa_x$ in \eqref{eq:kappadef} is a random variable, the autocorrelation is
\begin{equation}
   A_x(\xi)=\langle a_x(\xi)\rangle=\frac{1}{2}\sum_{\mu=1}^\infty\langle\hat{\kappa}_{x,2\mu}^2\rangle
        \cos(2\mu\pi\xi/{L}),\quad \xi\in[0, L].
\label{eq:autocorrx}
\end{equation}
The mean values of $\hat{\kappa}_{x,2\mu}^2$ are the coefficients in the expansion of
$A_x(\xi)$ in a cosine series such that
\begin{equation}
   \langle\hat{\kappa}_{x,2\mu}^2\rangle=\frac{4}{L}\int_0^L A_x(\xi)\cos(2\mu\pi\xi/{L})\, \ud \xi \equiv \hat{A}_{x\mu}
\label{eq:autocorrx2}
\end{equation}
as in Wiener-Khinchin's theorem.

As an example take
\begin{equation}
      A_x(\xi)=\left\{\begin{array}{rl}
                                      (1-\xi/\alpha)/\alpha,& \xi\in[0, \alpha],\\
                                      (1+\xi/\alpha)/\alpha,& \xi\in[-\alpha, 0),\\
                                      0,& {\rm otherwise},
                      \end{array}\right.
\label{eq:Axex}
\end{equation}
with $\alpha\in (0, L)$ which is $L$-periodic $A_x(\xi+L)=A_x(\xi)$. When $\alpha$ is small there is a correlation between the perturbations only in the vicinity.
The function in \eqref{eq:Axex} is scaled such that $\int_0^L A_x(\xi)\, d\xi =\int_{-L/2}^{L/2} A_x(\xi)\, d\xi =1$. The Fourier coefficients are
\begin{equation}
    \hat{A}_{x\mu}=2\cdot\frac{4}{L}\int_0^\alpha A_x(\xi)\cos(2\mu\pi\xi/{L})\, d\xi
                 =\frac{4L}{(\mu\pi\alpha)^2}\sin^2(\mu\pi\alpha/{L})> 0.
\label{eq:alphaex}
\end{equation}
The coefficients decay as $\mu^{-2}$ for increasing $\mu$ and when $\alpha$ is small then $A_x(\xi)$ approaches the
Dirac measure and $\hat{A}_{x\mu}=4/L+\ordo(\alpha^2)$.

There is a shift in the random frequency $\eta$ in \eqref{eq:deltadef} caused by $\delta_x$.
With a $\delta_x$ independent of $x$, the autocorrelation of $\delta_x$ is
\begin{equation}
     A_x(\xi)=\langle \delta_x(x+\xi)\delta_x(x)\rangle=\langle\delta_x^2\rangle=\veps^2\hG_{11I}^2\langle\tilde{\kappa}_{x1}^2\rangle
             =\frac{1}{4}\veps^2\hG_{11I}^2\langle\hat{\kappa}_{x2}^2\rangle
             =\frac{1}{4}\veps^2\hG_{11I}^2\hat{A}_{x1},
\label{eq:autocorrx3}
\end{equation}
and the deviation $s_x$ of the frequency from $\theta_1$ is
\begin{equation}
    s_x=\sqrt{A_x(0)}=\frac{1}{2}\veps\hG_{11I}\sqrt{\hat{A}_{x1}}.
\label{eq:apprdevfreqx}
\end{equation}



\section{Comparison between analysis and simulations}\label{sec:comparisons}
The PDE in \eqref{eq:Huang_model}, which is generalised to allow fluctuating coefficients $\fatsigma(x, t)$
depending on the realisation, is solved numerically in 1D in $[0, L]$.
The space derivative is approximated on a grid $x_i,\, i=0,1,\ldots,N,$ with constant grid size $\Delta x=x_{i}-x_{i-1}=L/N$ and
the usual difference formula of second order accuracy.
All reaction coefficients $\fatsigma$ are perturbed in time or in space and by the same factor such that $\fatf=0$ in \eqref{eq:Huang_1D},
$J=G$, and $H_1=G-(\gamma\pi^2/L^2)D$ in \eqref{eq:drhosol}.
The system in \eqref{eq:Huang_model} is solved by a Runge-Kutta method of fourth order accuracy with a constant time step and $N=21$.
This is a slight abuse of a numerical method designed for an ordinary differential equation. 
Fourth order temporal accuracy will not be achieved since $Y(t)$ from the OU process in Section~\ref{sec:OU} is only continuous.
An example of a solution is found in  Figure~\ref{fig:mincdepde}(a) with unperturbed coefficients.
The peak of $\rho_d$ alternates regularly between $x=0$ and $x=L$ with the frequency $\theta_1\approx 2\pi/40$.
The amplitude is approximately constant after an initial transient. The computed eigenvalues of $H_1$ and $H_{\omega}$ satisfy the assumption in \eqref{eq:lamassump}. A typical
cell cycle between the cell divisions of an {\it E. coli} is about 1200 s.


The frequency of the oscillations
is computed as in \eqref{eq:deltadefj} in an interval $[t_{j-1},\, t_{j+1}]$ of length $\Delta t_j$
where $t_j,\, j=0,1,\ldots ,$ are the time points of the consecutive extrema (maxima and minima) of the computed oscillations.
By \eqref{eq:deltadef} (and \eqref{eq:deltadefj} where $t_j\rightarrow t_{j+1}$) we have
\begin{equation}
     \eta_j=\theta_1 +\frac{1}{2}\veps\hG_{11I}\hat{\kappa}_{x2}+\veps\hG_{11I}\frac{1}{\Delta t_j}\int_{t_{j-1}}^{t_{j+1}} \kappa_t(v)\, \ud v
        = \theta_1(1+\veps\frac{\hG_{11I}}{\theta_1}\left(\frac{1}{2}\hat{\kappa}_{x2}+\bar{\kappa}_{tj})\right)).
\label{eq:avetaj}
\end{equation}
For the coefficients in \eqref{eq:Huang_param}, $\hG_{11I}\approx \theta_1$. Consequently,
\begin{equation}
     \eta_j\approx \theta_1(1+\veps\left(\frac{1}{2}\hat{\kappa}_{x2}+\bar{\kappa}_{tj}\right)).
\label{eq:avetaj2}
\end{equation}
The relative change in the frequency of the solution is $\veps(\frac{1}{2}\hat{\kappa}_{x2}+\bar{\kappa}_{tj})$.

The oscillatory frequency of the numerical $\rho_{d}$ solution due to temporal perturbations is compared at $x=0$ with the analysis in Section~\ref{sec:perturb}.
The OU process  (\eqref{eq:OU}) generates $\kappa_t(t)=Y(t)-1$ following Section~\ref{sec:OU}.
The diffusion $c$ for different relaxation parameters $\tau$ is chosen such that
\[
\frac{c\tau}{2}=\log 2.
\]
\begin{figure}[htbp]
\centering
\subfigure[$\tau=1$]{\includegraphics[width=0.7\linewidth]{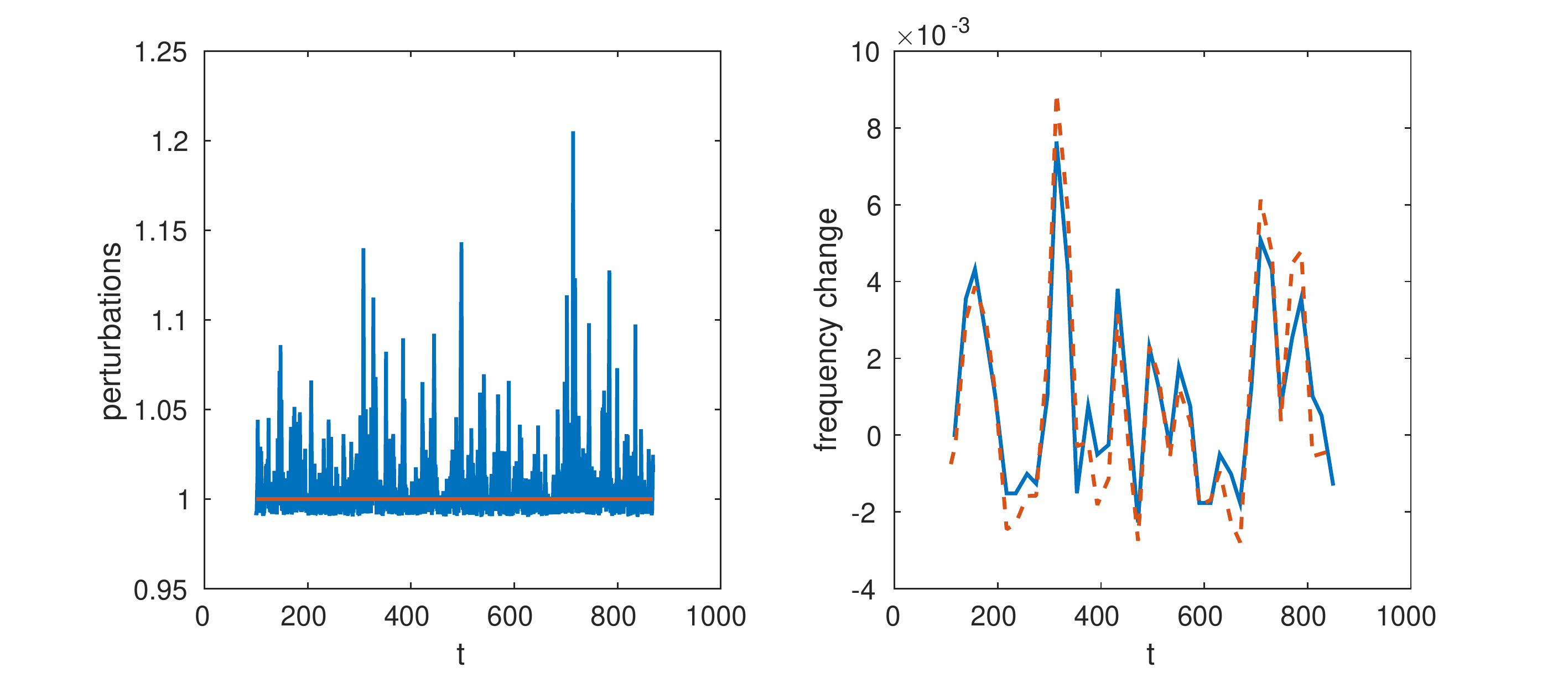}}
\subfigure[$\tau=10$]{\includegraphics[width=0.7\linewidth]{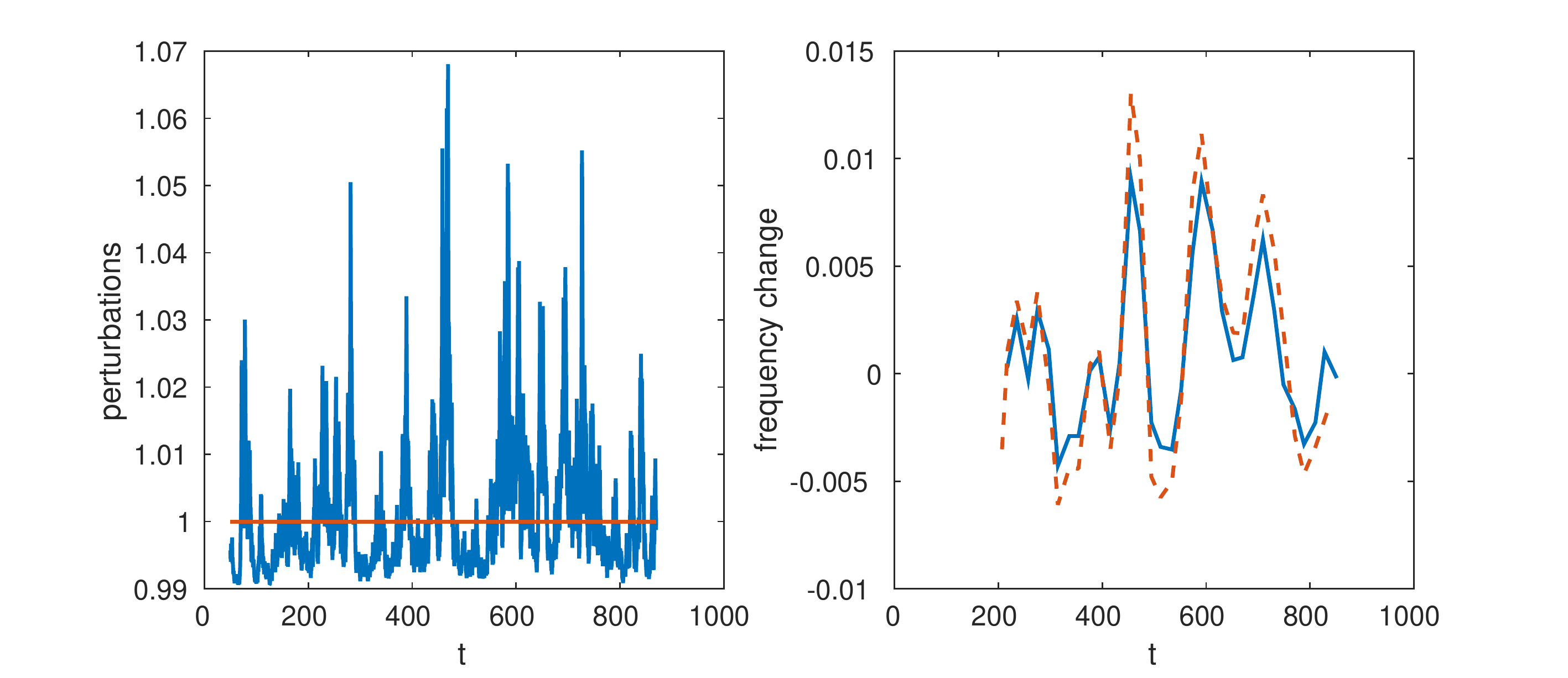}}
\subfigure[$\tau=100$]{\includegraphics[width=0.7\linewidth]{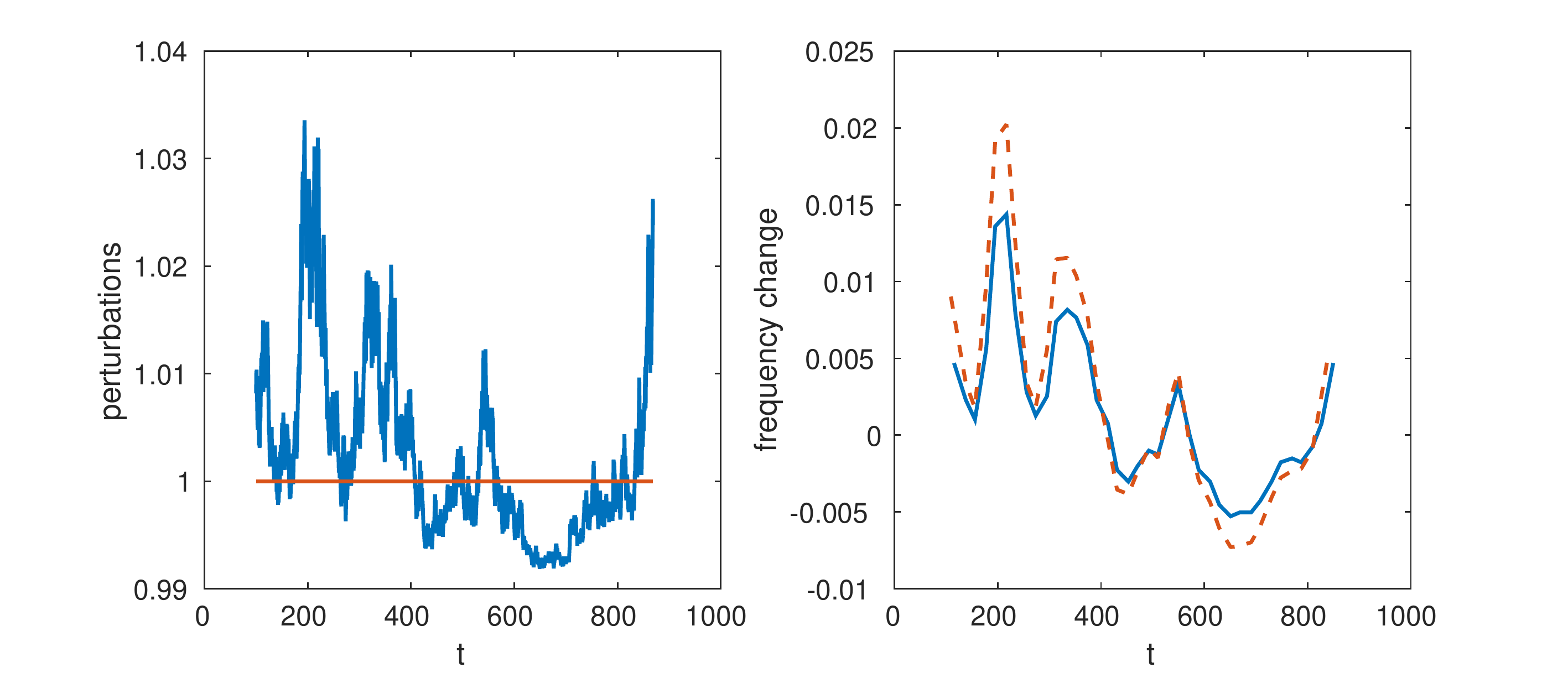}}
\subfigure[$\tau=1000$]{\includegraphics[width=0.7\linewidth]{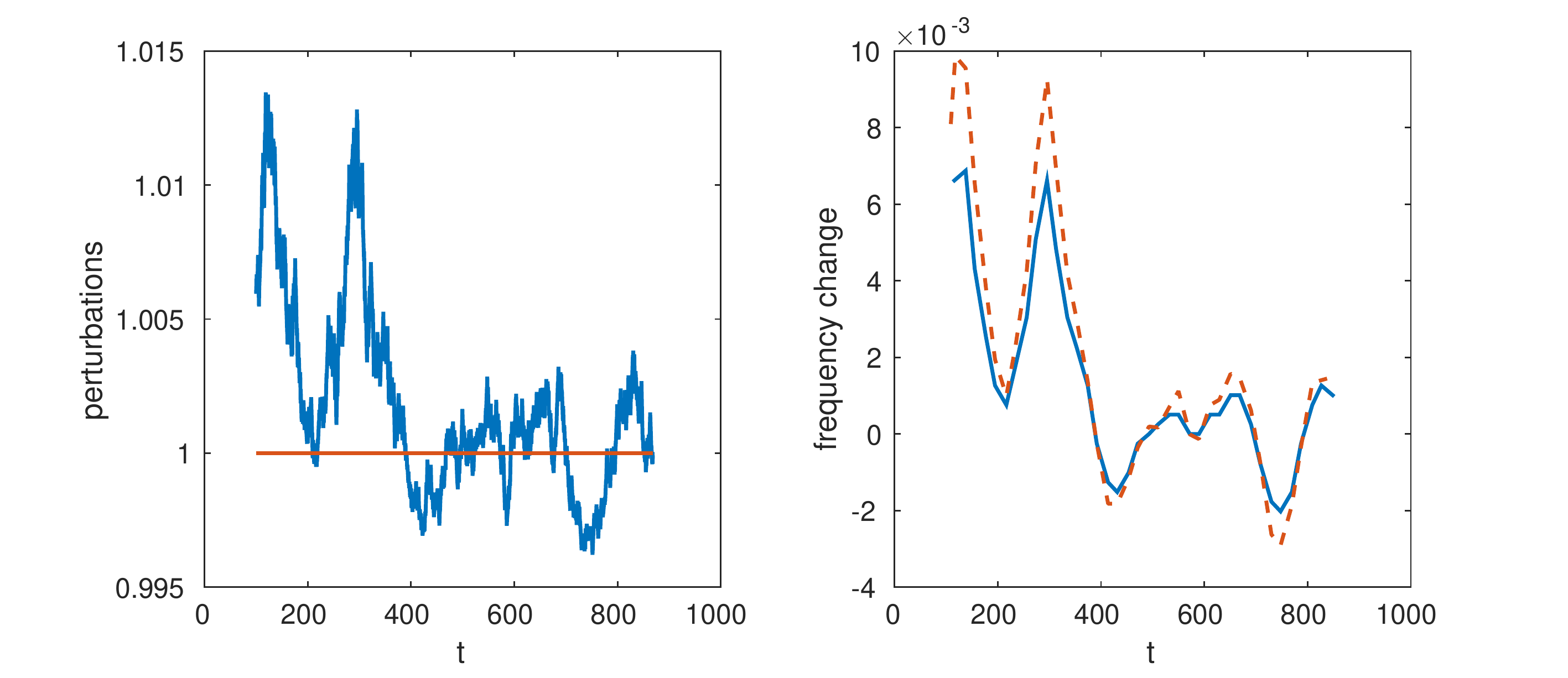}}
\caption{The temporal perturbations of the reaction coefficients are generated with $\veps=0.01$ by different $\tau$ parameters in the OU process. 
Left column: The  value of $1+\veps(Y(t)-1)$ generated by the OU process.
Right column: The relative change in the instantaneous frequency $\bar{\delta}_t(t)$ in \eqref{eq:deltatapprox} (solid blue)
and $\bar{\delta}_{OU}(t)$ in \eqref{eq:OUeffect} (dashed red). }
\label{fig:temp_perturb}
\end{figure}
Since
\[
\eta_j\Delta t_j=\theta_1(1+\delta_t/\theta_1)\Delta t_j=2\pi
\]
and $\theta_1T=2\pi$, the relative perturbation in the computed frequency in the nonlinear equations in $[t_{j-1}, t_{j+1}]$ is approximated as
\begin{equation}
     \frac{\delta_t(t_j)}{\theta_1}\approx \frac{\bar{\delta}_t(t_j)}{\theta_1}\equiv\frac{2\pi}{\theta_1 \Delta t_{j}}-1=\frac{T-\Delta t_j}{\Delta t_j},\, j\ge 1.
\label{eq:deltatapprox}
\end{equation}
The frequency $\theta_1$ is given by the unperturbed oscillations between two maxima or two minima.
The effect of the OU perturbations on the frequency in the same interval is
\begin{equation}
     \bar{\delta}_{OU}(t_j)=\veps\frac{\theta_1}{\Delta t_{j}}\int_{t_{j-1}}^{t_{j+1}} Y(v)-1 \, \ud v,
\label{eq:OUeffect}
\end{equation}
according to the analysis and \eqref{eq:avetaj2}. The quantities $\bar{\delta}_t(t_j)/\theta_1$ and $\bar{\delta}_{OU}(t_j)/\theta_1$ are
compared in the right column of Figure~\ref{fig:temp_perturb}
with good agreement.

\begin{figure}[htbp]
\includegraphics[width=0.33\linewidth]{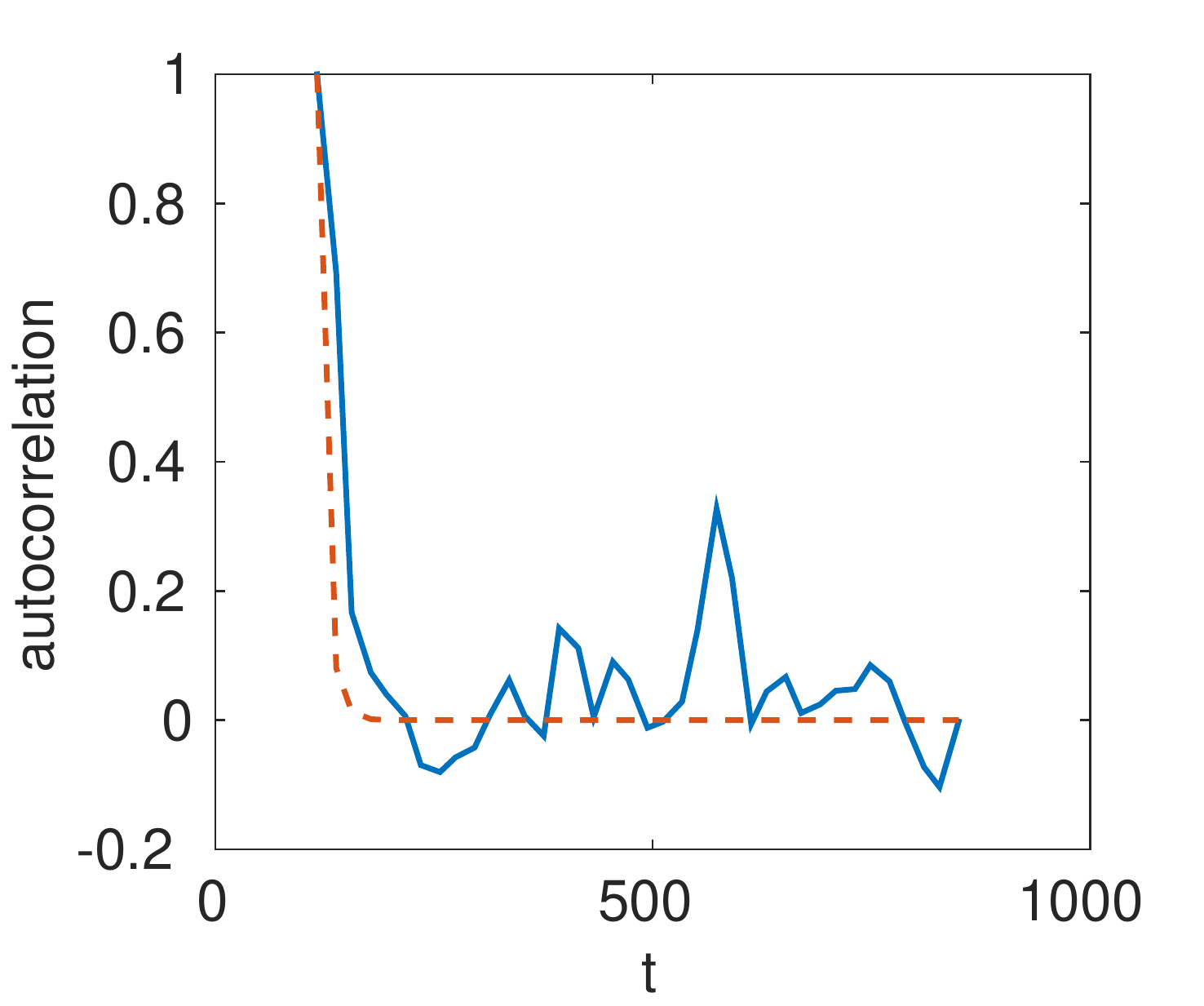}
\includegraphics[width=0.33\linewidth]{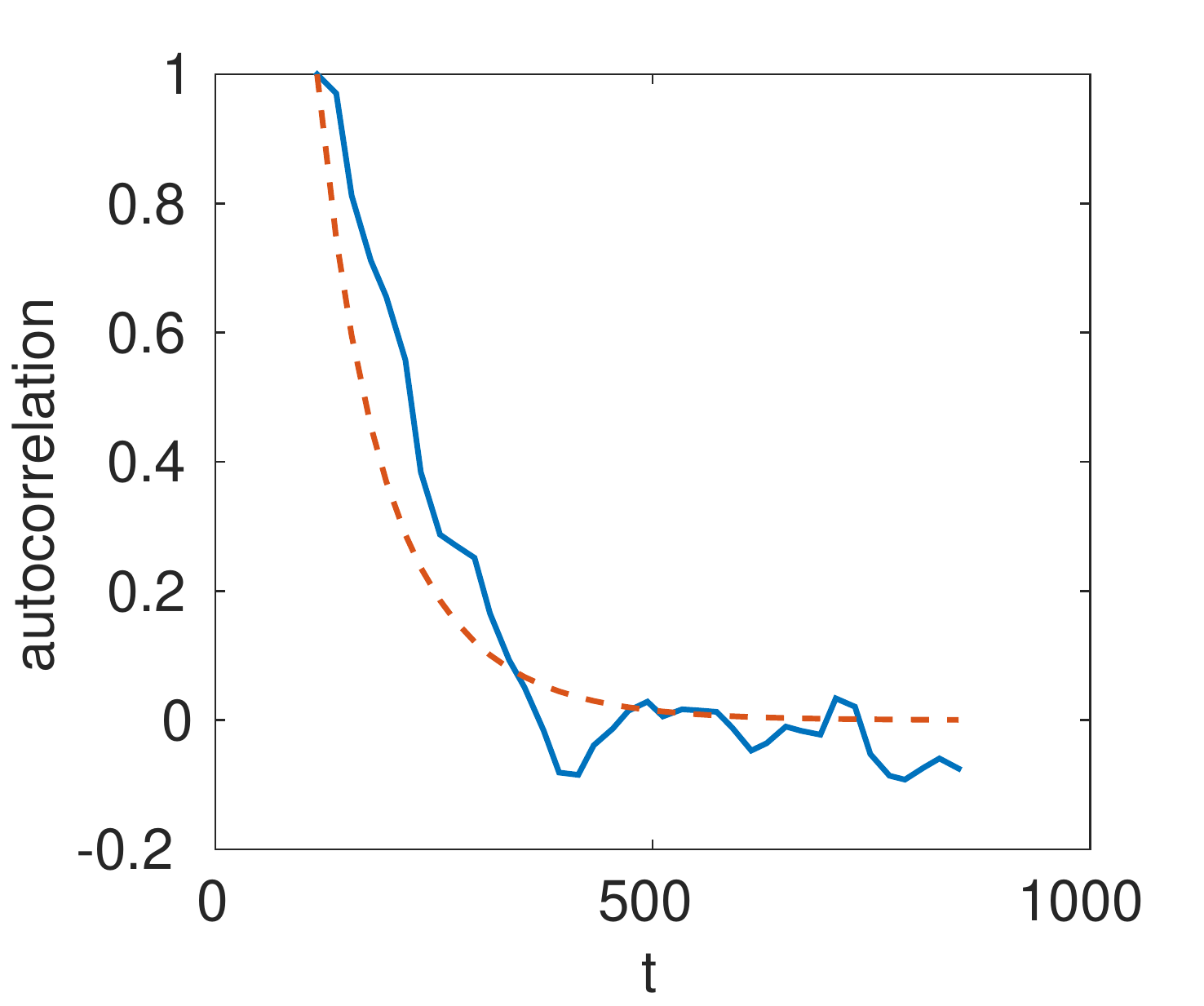}
\includegraphics[width=0.33\linewidth]{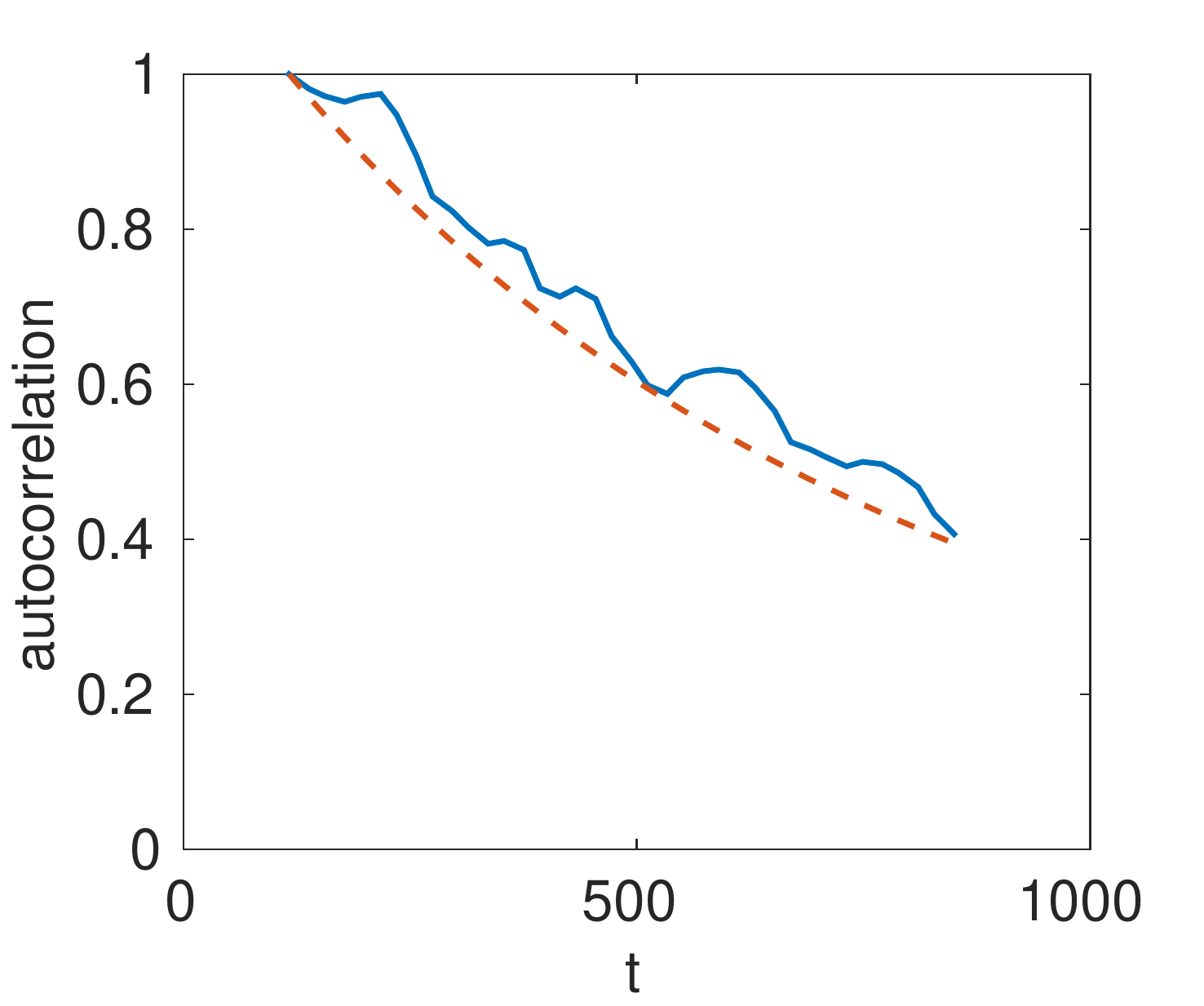}
\caption{The scaled autocorrelation ${A}_t(t_1,t)$ of the perturbations in the frequency in \eqref{eq:autocorrtj} in the PDE solution (solid blue)
for different $\tau$ parameters in the OU process with $\veps=0.01$ compared with the estimate in \eqref{eq:autocorrtj} (dashed red).
Left: $\tau=10$. Middle: $\tau=100$.
Right: $\tau=1000$. }
\label{fig:temp_perturb_auto}
\end{figure}

The autocorrelation of the observed instantaneous frequency change in \eqref{eq:autocorrtj}
${A}_t(t_1, t_k)\approx \veps^2\theta_1^2\langle \kappa_t(t_1)\kappa_t(t_k)\rangle$ is compared to the estimate in \eqref{eq:autocorrtj}
in Figure~\ref{fig:temp_perturb_auto} for three different relaxations $\tau$. The $\langle\; \cdot\;\rangle$ average is taken over 200 trajectories and
the data are scaled by the initial ${A}_t(t_1, t_1)$.
Since $\exp(c\tau/2)=2$, the correlation depends only on the time scale $\tau$.
The estimate behaves as $1-\frac{\log 2}{2}\frac{\Delta t_{1k}}{\tau}$ for small $\Delta t_{1k}/\tau$ and is very small for large $\Delta t_{1k}/\tau$,
see \eqref{eq:autocorrtappr1} and \eqref{eq:autocorrtappr2}. The transient phase is short for $\tau=1$ and not over at 1000 s when $\tau=1000$.

\begin{figure}[htbp]
\begin{center}
\includegraphics[width=0.8\linewidth]{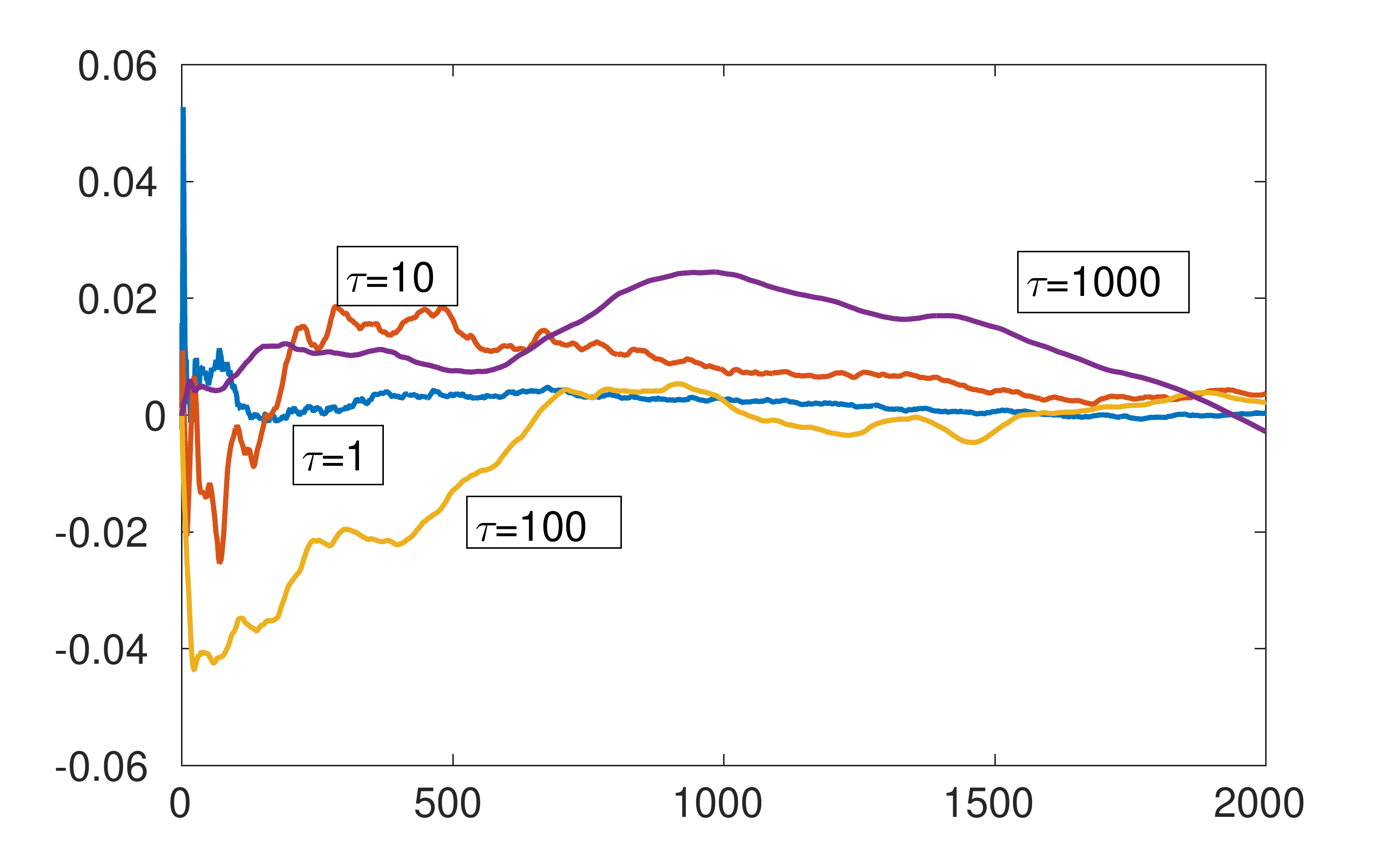} \\
\includegraphics[width=0.8\linewidth]{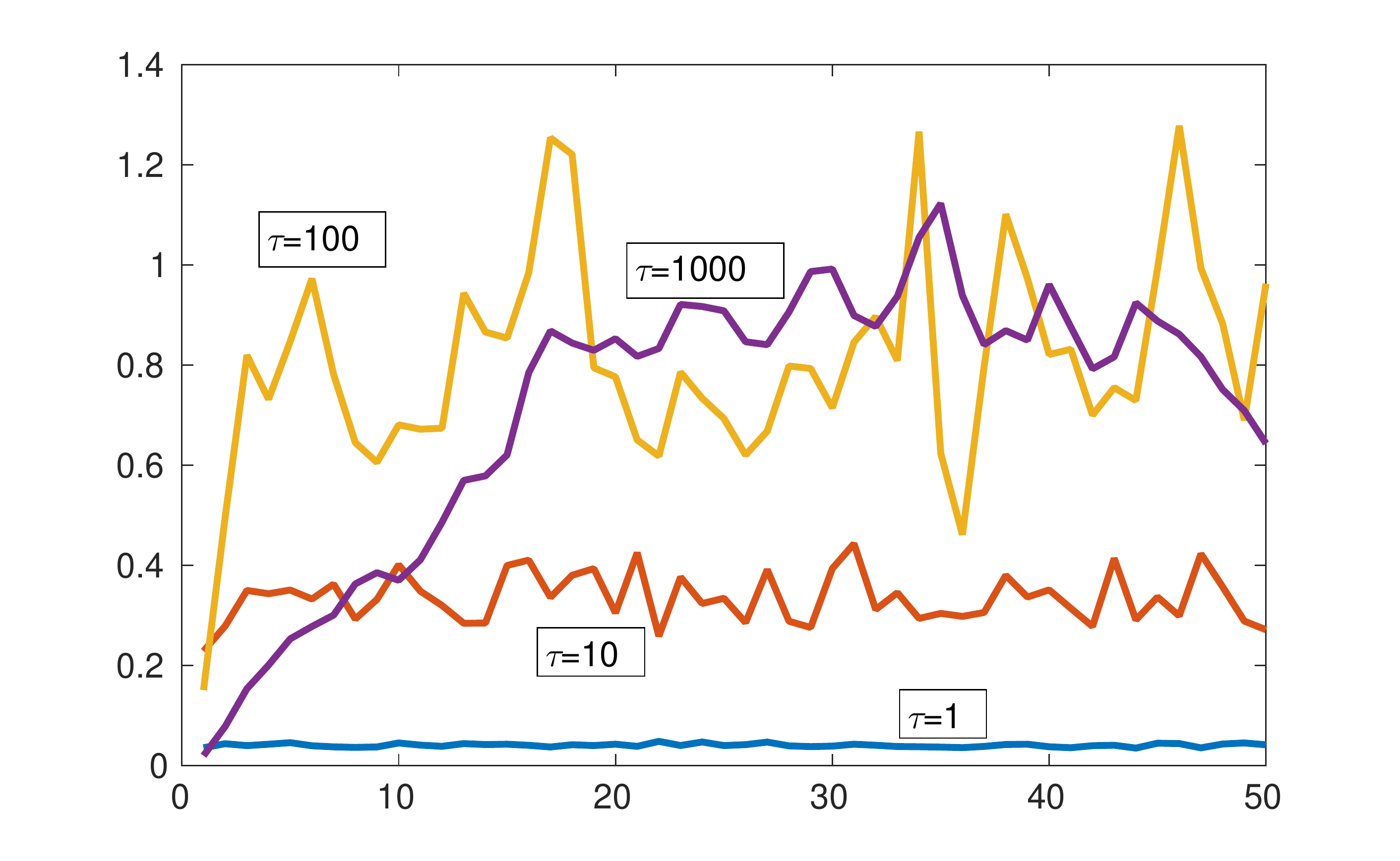} \\
\includegraphics[width=0.8\linewidth]{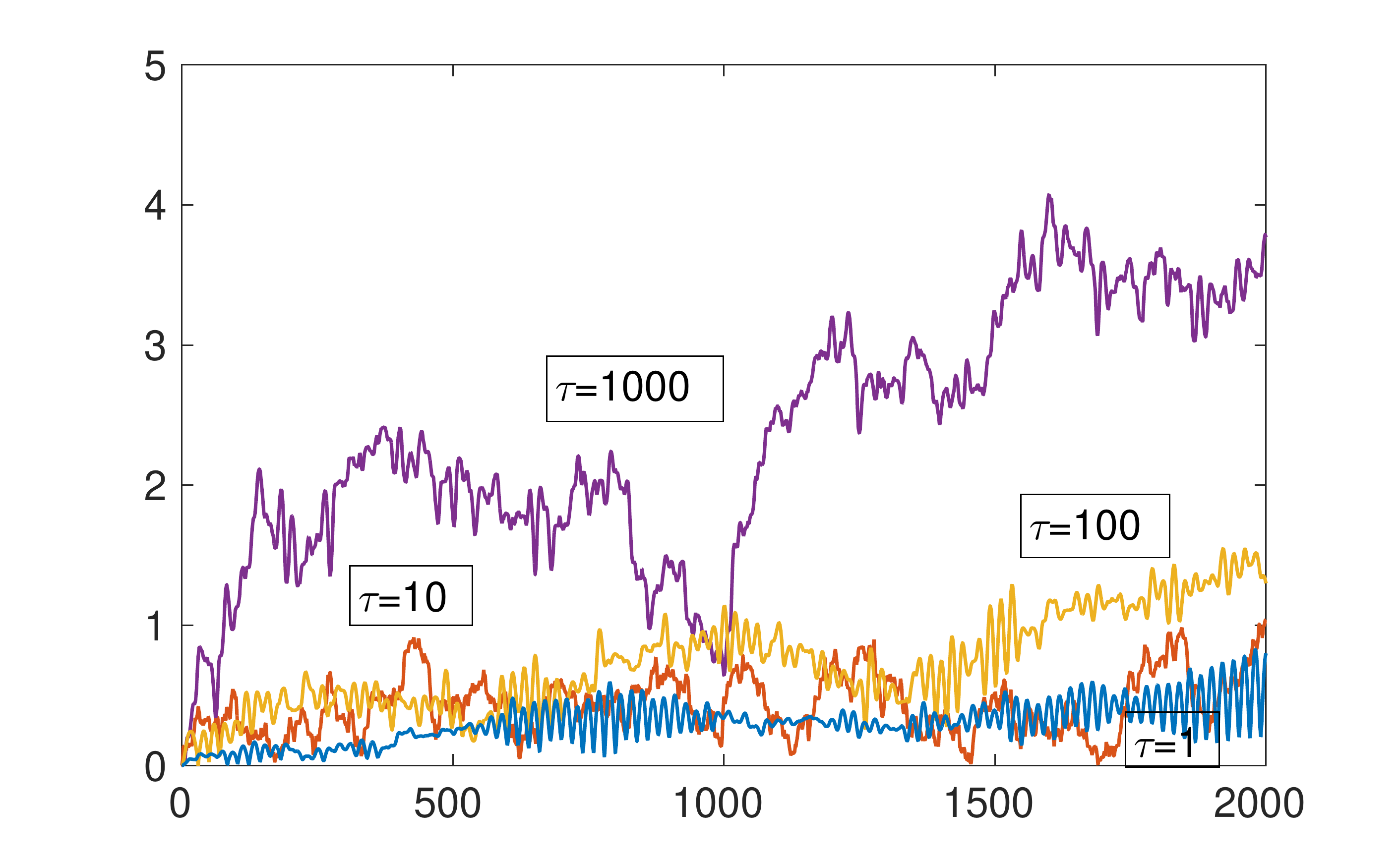}
\end{center}
\caption{Averaged integrals of $\kappa_t(t)$ over 400 realizations for different $\tau$ ($\tau=1$: blue; $\tau=10$: red; $\tau=100$: yellow; $\tau=1000$: purple).
Top: $\frac{1}{t}\int_0^t \langle\kappa_t(s)\rangle\, \ud s$ for $t\in [0, 2000]$.
Middle:           $\langle\bar{\kappa}^2_{tj}\rangle=\langle\left(\frac{1}{\Delta t_j}\int^{t_{j+1}}_{t_{j-1}} \kappa_t(s)\, \ud s\right)^2\rangle$ in 50 intervals
of length about 40 covering an interval of length about 2000.
Bottom: $|\exp(i\theta_1 t)\int_0^t \exp(-2i\theta_1 s)\langle\kappa_t(s)\rangle \,\ud s|$  for $t\in [0, 2000]$.}
\label{fig:kappat}
\end{figure}

The properties of the temporal perturbation $\kappa_t$ generated by the OU process are evaluated in Figure~\ref{fig:kappat} for different $\tau$. Averages of integrals
of $\kappa_t$ are taken over 400 trajectories. After a transient phase, a stationary distribution $p_\infty$ of $\kappa_t$ is obtained, cf. Section~\ref{sec:OU}.
The transient is longer the larger $\tau$ is as in Figure~\ref{fig:temp_perturb_auto}.
In the upper left figure, the integral tends to 0 as required in \eqref{eq:enoise} but the convergence is fast when
$\tau=1$ and slow for $\tau=1000$. The response of the OU perturbations in the frequency $\eta(t)$ in \eqref{eq:deltadef} vanishes for large $t$ but at different speed
depending on $\tau$.
The average of the square of $\kappa_t$ over a period in \eqref{eq:avetaj}
in the upper right figure is small for $\tau=1$ and growing for increasing $\tau$.
This is also the trend in Figure~\ref{fig:temp_perturb} for the amplitude of the frequency change.
The average $\langle\kappa_t^2\rangle$ estimates the variance of $\kappa_t$ after the transient phase since $\langle\bar{\kappa}_t\rangle\approx 0$ there.
When $\tau=1$ the fast fluctuations in $\kappa_t(t)$,
see e.g. the left column of Figure~\ref{fig:temp_perturb}, are averaged efficiently over $T$ in $\bar{\kappa}_t$.
The interval is not sufficiently long for evaluation of the average for $\tau=1000$. The absolute value of the integral in the lower panel causes part of the change in the
amplitude of the oscillations in \eqref{eq:du1sol}. Also here there is an effect of the $\tau$ parameter in the OU process.


The shift in the frequency caused by the spatial perturbations in \eqref{eq:deltadef} and \eqref{eq:avetaj2} is $\delta_x=\frac{1}{2}\veps\theta_1\hat{\kappa}_{x2}$.
In the first experiment, the perturbations of the reaction coefficients are smooth with
$1+\veps\kappa_x(x)=1+\veps\cos(\omega\pi x/L)$ where $\veps=0.02$ and $0.1$
and $\omega=2$. The measured change in frequency is $\bar{\delta}_x$ in the solution of \eqref{eq:Huang_model} and is computed as in \eqref{eq:deltatapprox}.
For $\veps=0.02$, $\bar{\delta}_x/\theta_1\approx 0.01$ and for $\veps=0.1$, $\bar{\delta}_x/\theta_1\approx 0.05$.
The expected relative frequency shift in \eqref{eq:deltatapprox} is $\delta_x/\theta_1=\frac{1}{2}\veps\hat{\kappa}_{x2}=\frac{1}{2}\veps$
which is in good agreement with the computed shifts.
With $\omega=6$ and $\veps=0.1$, we have $\bar{\delta}_x/\theta_1=0.011$ in a numerical experiment which is $1/5$ of the value at $\omega=2$.

\begin{figure}[htbp]
\begin{center}
\includegraphics[width=0.45\linewidth]{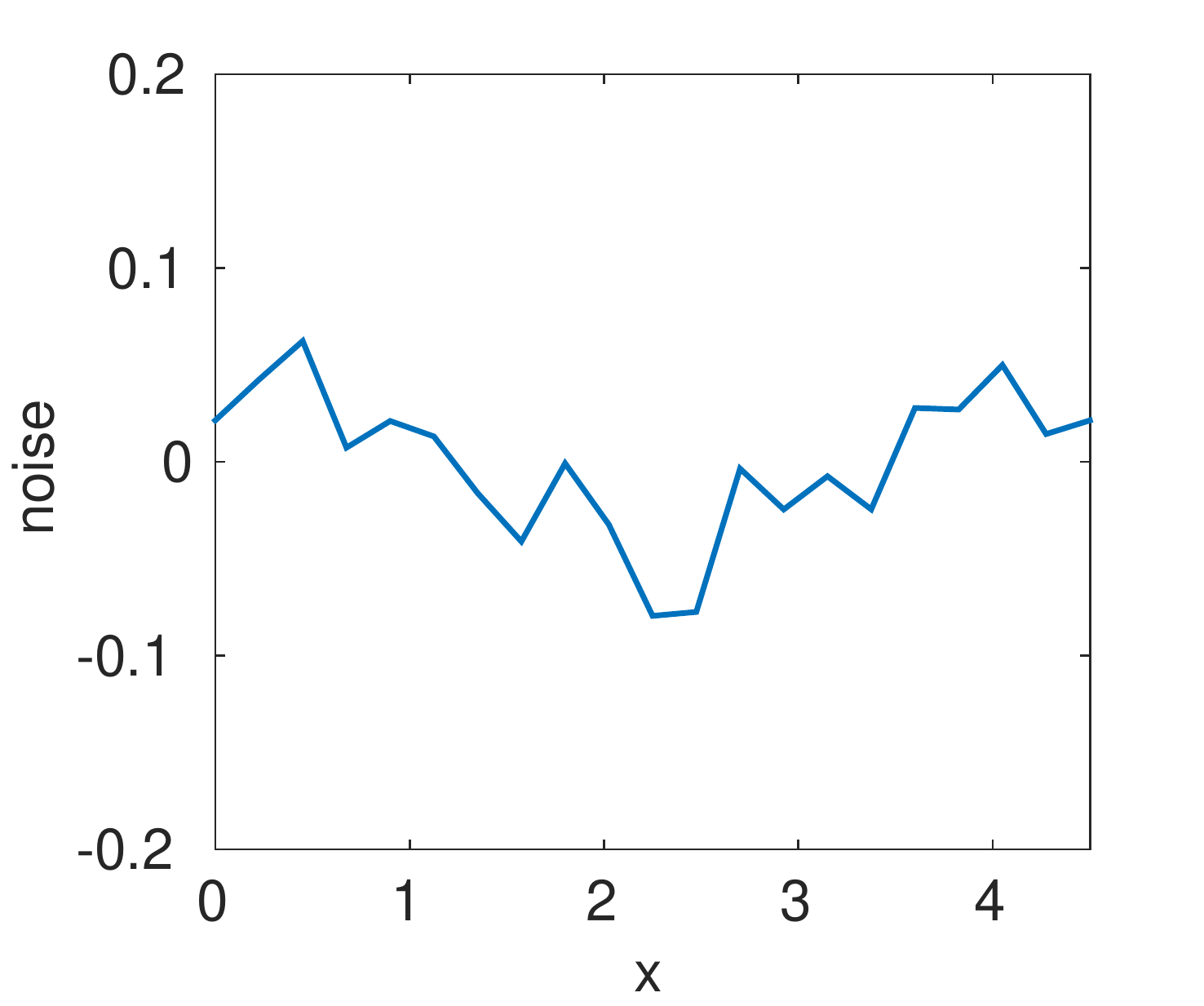}
\includegraphics[width=0.45\linewidth]{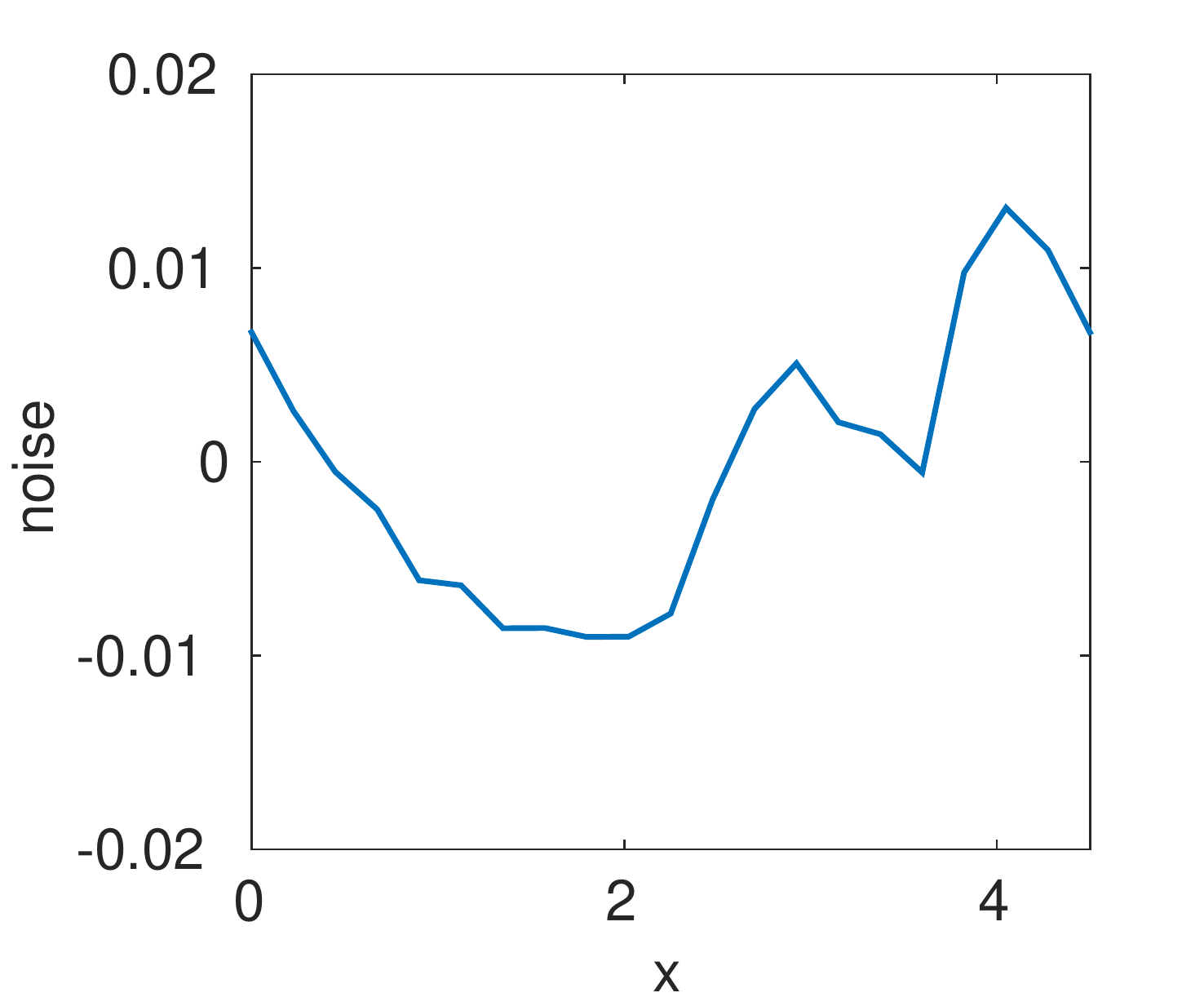}
\\
\includegraphics[width=0.45\linewidth]{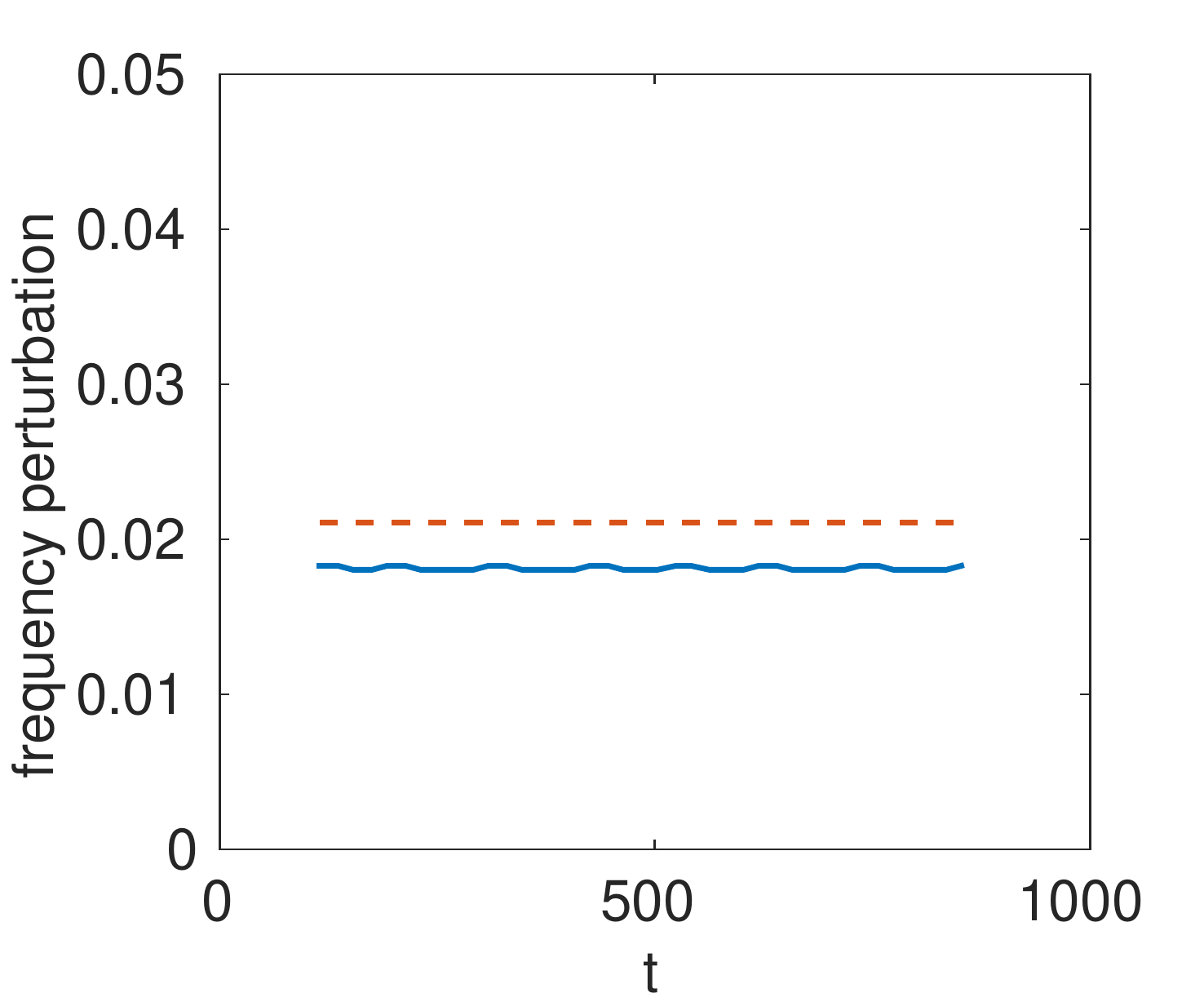}\hspace{0.3cm}
\includegraphics[width=0.45\linewidth,angle=0,origin=c]{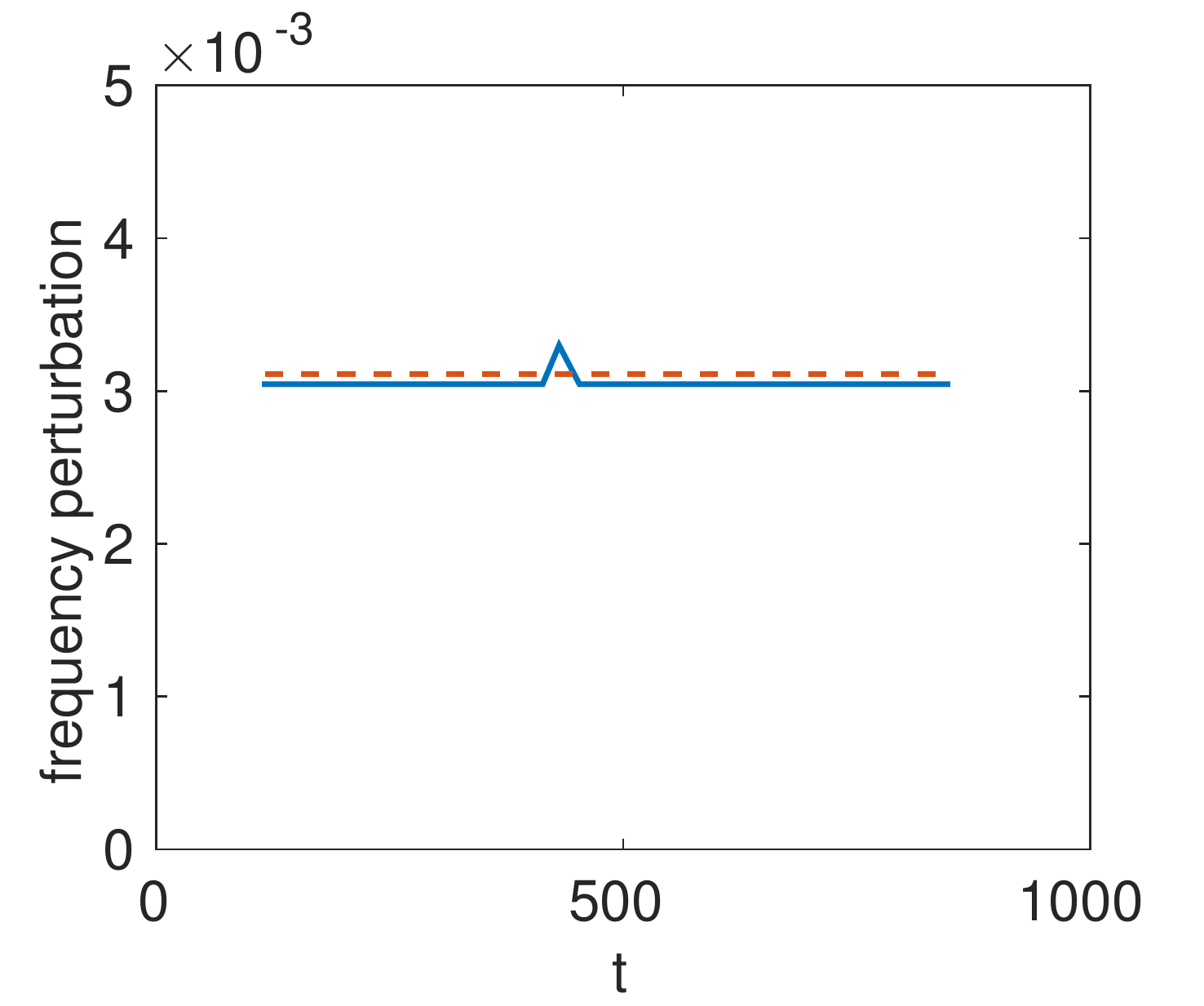} 
\end{center}
\caption{The spatial perturbations of the reaction coefficients are generated for $x\in [0, 4.5]$ by a multivariate normal distribution with $\alpha=2$ in \eqref{eq:Axex}.
Left column: $\veps=0.1$.
Right column: $\veps=0.01$.
Upper row: The $\epsilon\kappa_x(x)$ value.
Lower row: The measured relative change $\delta_x/\theta_1$ (solid blue) in the instantaneous frequency in the time interval and the prediction by linear theory (dashed red). }
\label{fig:spatial_perturb_rand}
\end{figure}

In the next experiment, the perturbations in the coefficients are random in space with zero mean.
The perturbations at $x_i$ are sampled from a multivariate normal distribution and they are correlated in neighboring grid points $x_i$ and $x_j$ as
described in Section~\ref{sec:spatnoise}. The correlations are a discretization of the ones in \eqref{eq:Axex} with a symmetric, circulant, and Toeplitz covariance matrix $C$.
For $\alpha<L$, the elements of $C$ are
\[
    C_{ij}=A_x(x_i-x_j)+A_x(x_i-x_j-L)+A_x(x_i-x_j+L),\; i,j=1,\ldots,N.
\]
The procedure in Section~\ref{sec:spatnoise} generates $\kappa_{xi},\, i=0,1,\ldots,N,$ for a given $\alpha$ in \eqref{eq:Axex}
and their mean $\overline{\kappa}_x$ is computed. Then $\kappa_x(x_i)=\kappa_{xi}-\overline{\kappa}_x$ such that $\kappa_x(x_i)$ has zero mean as in \eqref{eq:enoise}.
The relative perturbation $\delta_x/\theta_1$ is compared to the predictions in \eqref{eq:avetaj2}.
The coefficient $\hat{\kappa}_{x2}$ is determined by the discrete cosine transform of $\kappa_x(x_i)$. In Figure~\ref{fig:spatial_perturb_rand},
two examples of perturbations are found with $\alpha=2$ and the corresponding change in frequency.
The frequency shifts are almost constant in time and the theory agrees well with the experiments.
Small oscillations are observed in the perturbed frequency due to the sensitivity to the computed time interval $\Delta t_j$ in \eqref{eq:deltatapprox}.
If the unperturbed interval is $\overline{\Delta t}=2\pi/\theta_1$ and the perturbed interval is $\Delta t_j=\overline{\Delta t}(1-\mu_j)$ with the
relative perturbation $\mu_j$, then $\bar{\delta}_t(t_j)$ in \eqref{eq:deltatapprox} is
\[
     \frac{\bar{\delta}_t(t_j)}{\theta_1}=\frac{2\pi}{\theta_1 \Delta t_{j}}-1=\frac{1}{1-\mu_j}-1\approx \mu_j.
\]
The relative precision in the numerical computations of $\Delta t_j$ has to be much better than $\mu_j$ which is about 0.02 and 0.003 in the figures.

The autocorrelation in \eqref{eq:autocorrx3} is determined for $\veps=0.1$ by averaging $(\delta_x/\theta_1)^2$ over 200 realizations resulting in $s_x/\theta_1=0.0203$.
The corresponding theoretical value of $\frac{1}{2}\veps\sqrt{\langle \hat{\kappa}_{x2}^2\rangle}$ is 0.0223 in a reasonable agreement between the nonlinear model and the
analysis of the linearization.


\section{Discussion}\label{sec:disc}

We have explored the robustness of spatiotemporal oscillations using a complementary combination of analysis and simulation.
All parameters in $\fatsigma$ in \eqref{eq:Huang_model} and \eqref{eq:Huang_param} are perturbed by $\veps\kappa_t(t)$ in time and $\veps\kappa_x(x)$ in space.
Numerical solutions of the nonlinear system of PDEs \eqref{eq:Huang_model} with extrinsic noise in $\fatsigma$ are compared to the solution without noise. The changes in
the oscillation frequency agree well with the theoretical predictions for a linearized system satisfying certain assumptions.
The MinD model in \eqref{eq:Huang_model} fulfills these assumptions and is robust in the sense that small perturbations 
in $\fatsigma$ result in small differences in the frequency and the amplitude. 

The analytical approach is suitable for many other systems satisfying the assumptions concerning the 
eigenvalues of the Jacobian of the linearized system \eqref{eq:Huang_eps} and that the parameters in the reaction rates should appear linearly in the equations. 
It is likely that the conclusions concerning the Ornstein-Uhlenbeck perturbations in time and the correlated spatial perturbations are
much more general than just for the particular example that we have studied here.



\section{Acknowledgment}
SM was supported by the Centre for Interdisciplinary Mathematics when visiting Uppsala and PL was supported by University of New South Wales during his visit there. AH was supported by the Swedish strategic research programme eSSENCE. 

\bibliographystyle{siamplain}
\bibliography{extrinsic}

\newcommand{\noopsort}[1]{}
\begin{thebibliography}{10}

\bibitem{BlaKae03}
{\sc W.~J. Blake, M.~Kaern, C.~R. Cantor, and J.~J. Collins}, {\em Noise in
  eukaryotic gene expression.}, Nature, 422 (2003), pp.~633--637,
  \href{http://dx.doi.org/10.1038/nature01546} {doi:10.1038/nature01546},
  \url{http://dx.doi.org/10.1038/nature01546}.

\bibitem{URDME}
{\sc B.~Drawert, S.~Engblom, and A.~Hellander}, {\em {URDME}: a modular
  framework for stochastic simulation of reaction-transport processes in
  complex geometries}, BMC Syst. Biol., 6 (2012), p.~76,
  \href{http://dx.doi.org/10.1186/1752-0509-6-76} {doi:10.1186/1752-0509-6-76}.

\bibitem{EloLev02}
{\sc M.~B. Elowitz, A.~J. Levine, E.~D. Siggia, and P.~S. Swain}, {\em
  Stochastic gene expression in a single cell.}, Science, 297 (2002),
  pp.~1183--1186, \href{http://dx.doi.org/10.1126/science.1070919}
  {doi:10.1126/science.1070919},
  \url{http://dx.doi.org/10.1126/science.1070919}.

\bibitem{FanElf06}
{\sc D.~Fange and J.~Elf}, {\em Noise-induced {M}in phenotypes in \emph{ {E.
  coli}}}, PLoS Comput. Biol., 2 (2006), pp.~0637--0647,
  \href{http://dx.doi.org/10.1371/journal.pcbi.0020080}
  {doi:10.1371/journal.pcbi.0020080}.

\bibitem{FedFon02}
{\sc N.~Fedoroff and W.~Fontana}, {\em Small numbers of big molecules},
  Science, 297 (2002), pp.~1129 -- 1131,
  \href{http://dx.doi.org/10.1126/science.1075988}
  {doi:10.1126/science.1075988}.

\bibitem{FiFrMeLuChKr2010}
{\sc E.~Fischer-Friedrich, G.~Meacci, J.~Lutkenhaus, H.~Chat{\'e}, and
  K.~Kruse}, {\em Intra- and intercellular fluctuations in {M}in-protein
  dynamics decrease with cell length}, Proc. Natl. Acad. Sci. USA, 107 (2010),
  pp.~6134--6139, \href{http://dx.doi.org/10.1073/pnas.0911708107}
  {doi:10.1073/pnas.0911708107}.

\bibitem{Gil92}
{\sc D.~Gillespie}, {\em Markov Processes: An Introduction for Physical
  Scientists}, Academic Press, 1992,
  \url{https://www.elsevier.com/books/markov-processes/gillespie/978-0-12-283955-9}.

\bibitem{gillespie}
{\sc D.~T. Gillespie}, {\em A general method for numerically simulating the
  stochastic time evolution of coupled chemical reactions}, J.~Comput.~Phys.,
  22 (1976), pp.~403--434,
  \href{http://dx.doi.org/10.1016/0021-9991(76)90041-3}
  {doi:10.1016/0021-9991(76)90041-3}.

\bibitem{Gil96}
{\sc D.~T. Gillespie}, {\em Exact numerical simulation of the
  {O}rnstein-{U}hlenbeck process and its integral}, Phys. Rev. E, 54 (1996),
  pp.~2084--2091, \href{http://dx.doi.org/10.1103/PhysRevE.54.2084}
  {doi:10.1103/PhysRevE.54.2084}.

\bibitem{Gould13}
{\sc P.~D. Gould, N.~Ugarte, M.~Domijan, M.~Costa, J.~Foreman, D.~MacGregor,
  K.~Rose, J.~Griffiths, A.~J. Millar, B.~Finkelst{\"a}dt, S.~Penfield, D.~A.
  Rand, K.~J. Halliday, and A.~J.~W. Hall}, {\em Network balance via {CRY}
  signalling controls the \emph{Arabidopsis} circadian clock over ambient
  temperatures}, Mol. Syst. Biol., 9 (2013), p.~650,
  \href{http://dx.doi.org/10.1038/msb.2013.7} {doi:10.1038/msb.2013.7}.

\bibitem{GunawardenaPatheitc}
{\sc J.~Gunawardena}, {\em Models in biology: `accurate descriptions of our
  pathetic thinking'}, BMC Biology, 12 (2014),
  \href{http://dx.doi.org/10.1186/1741-7007-12-29}
  {doi:10.1186/1741-7007-12-29}.

\bibitem{Halatek12}
{\sc J.~Halatek and E.~Frey}, {\em Highly canalized {MinD} transfer and {MinE}
  sequestration explain the origin of robust {MinCDE}-protein dynamics}, Cell,
  1 (2012), pp.~741--52, \href{http://dx.doi.org/10.1016/j.celrep.2012.04.005}
  {doi:10.1016/j.celrep.2012.04.005}.

\bibitem{HilfingerPaulsson2011}
{\sc A.~Hilfinger and J.~Paulsson}, {\em Separating intrinsic from extrinsic
  fluctuations in dynamic biological systems}, Proc. Acad. Natl. Sci., 109
  (2011), pp.~12167--72, \href{http://dx.doi.org/10.1073/pnas.1018832108}
  {doi:10.1073/pnas.1018832108}.

\bibitem{Huang:2003}
{\sc K.~C. Huang, Y.~Meir, and N.~S. Wingreen}, {\em Dynamic structures in
  \emph{Escherichia coli}: {S}pontaneous formation of {M}in{E} and {M}in{D}
  polar zones}, Proc.~Natl.~Acad.~Sci.~USA, 100 (2003), pp.~12724--12728,
  \href{http://dx.doi.org/10.1073/pnas.2135445100}
  {doi:10.1073/pnas.2135445100}.

\bibitem{Kerr06}
{\sc R.~A. Kerr, H.~Levine, T.~J. Sejnowski, and W.-J. Rappel}, {\em Division
  accuracy in a stochastic model of {M}in oscillations in \emph{Escherichia
  coli}}, Proc. Natl. Acad. Sci. USA, 103 (2006), pp.~347--352,
  \href{http://dx.doi.org/10.1073/pnas.0505825102}
  {doi:10.1073/pnas.0505825102}.

\bibitem{Kirschner50YearsJacobMonod2011}
{\sc M.~Kirschner, L.~Shapiro, H.~McAdams, G.~Almouzni, P.~Sharp, R.~Young, and
  U.~Alon}, {\em Fifty years after {J}acob and {M}onod: what are the unanswered
  questions in molecular biology?}, Mol. Cell, 42 (2011), pp.~403--4,
  \href{http://dx.doi.org/http://dx.doi.org/10.1016/j.molcel.2011.05.003}
  {doi:http://dx.doi.org/10.1016/j.molcel.2011.05.003}.

\bibitem{Kruse:2002}
{\sc K.~Kruse}, {\em A dynamic model for determining the middle of
  \emph{Escherichia coli}}, Biophys.~J., 82 (2002), pp.~618--627,
  \href{http://dx.doi.org/10.1016/S0006-3495(02)75426-X}
  {doi:10.1016/S0006-3495(02)75426-X}.

\bibitem{KruHowMar:2007}
{\sc K.~Kruse, M.~Howard, and W.~Margolin}, {\em An experimentalist's guide to
  computational modelling of the {M}in system}, Mol.~Microb., 63 (2007),
  pp.~1279--1284, \href{http://dx.doi.org/10.1111/j.1365-2958.2007.05607.x}
  {doi:10.1111/j.1365-2958.2007.05607.x}.

\bibitem{McAArk97}
{\sc H.~H. McAdams and A.~Arkin}, {\em Stochastic mechanisms in gene
  expression}, Proc.~Natl.~Acad.~Sci.~USA, 94 (1997), pp.~814--819,
  \url{http://www.pnas.org/content/94/3/814.abstract}.

\bibitem{MeaKru:2005}
{\sc G.~Meacci and K.~Kruse}, {\em Min oscillations in \emph{Escherichia coli}
  induced by interactions of membrane-bound proteins}, Phys.~Biol., 2 (2005),
  pp.~89--97, \href{http://dx.doi.org/10.1088/1478-3975/2/2/002}
  {doi:10.1088/1478-3975/2/2/002}.

\bibitem{Mur02}
{\sc J.~Murray}, {\em Mathematical biology : an introduction}, New York :
  Springer, 2002, \url{http://www.springer.com/gp/book/9780387952239}.

\bibitem{Nayfeh}
{\sc A.~H. Nayfeh}, {\em Perturbation Methods}, Wiley, 1973,
  \href{http://dx.doi.org/10.1002/9783527617609} {doi:10.1002/9783527617609}.

\bibitem{Oks98}
{\sc B.~K. {\O}ksendal}, {\em Stochastic Differential Equations: An
  Introduction with Applications}, Berlin ; New York : Springer, 1998,
  \href{http://dx.doi.org/10.1007/978-3-642-14394-6}
  {doi:10.1007/978-3-642-14394-6}.

\bibitem{PhilipsPhysBiologyCellBook2012}
{\sc R.~Phillips, J.~Kondev, and J.~Theriot}, {\em Physical Biology of the
  Cell}, Garland Science, 2012,
  \url{http://www.garlandscience.com/product/isbn/9780815344506}.

\bibitem{RenRuo02}
{\sc L.~Rensing and P.~Ruoff}, {\em Temperature effect on entrainment, phase
  shifting, and amplitude of circadian clocks and its molecular bases},
  Chronobiol. Int., 19 (2002), pp.~807--864,
  \href{http://dx.doi.org/10.1081/CBI-120014569} {doi:10.1081/CBI-120014569}.

\bibitem{ShaOll08}
{\sc V.~Shahrezaei, J.~Ollivier, and P.~Swain}, {\em Colored extrinsic
  fluctuations and stochastic gene expression}, Mol. Syst. Biol., 4 (2008),
  pp.~1--9, \href{http://dx.doi.org/10.1038/msb.2008.31}
  {doi:10.1038/msb.2008.31}.

\bibitem{ShaSwa08}
{\sc V.~Shahrezaei and P.~S. Swain}, {\em The stochastic nature of biochemical
  networks}, Curr. Op. Biotech., 19 (2008), pp.~369--374,
  \href{http://dx.doi.org/10.1016/j.copbio.2008.06.011}
  {doi:10.1016/j.copbio.2008.06.011}.

\bibitem{SteGilDoy2004}
{\sc J.~Stelling, E.~D. Gilles, and F.~D. III}, {\em Robustness properties of
  circadian clock architectures}, Proc. Natl. Acad. Sci. USA, 101 (2004),
  pp.~13210--13215, \href{http://dx.doi.org/10.1073/pnas.0401463101}
  {doi:10.1073/pnas.0401463101}.

\bibitem{StellingDoyleRobustness2004}
{\sc J.~Stelling, U.~Sauer, Z.~Szallasi, F.~J.~D. III, and J.~Doyle}, {\em
  Robustness of cellular functions}, Cell, 118 (2004), pp.~675--685,
  \href{http://dx.doi.org/10.1016/j.cell.2004.09.008}
  {doi:10.1016/j.cell.2004.09.008}.

\bibitem{SwaElo02}
{\sc P.~S. Swain, M.~B. Elowitz, and E.~D. Siggia}, {\em Intrinsic and
  extrinsic contributions to stochasticity in gene expression.},
  Proc.~Natl.~Acad.~Sci.~USA, 99 (2002), pp.~12795--12800,
  \href{http://dx.doi.org/10.1073/pnas.162041399} {doi:10.1073/pnas.162041399},
  \url{http://dx.doi.org/10.1073/pnas.162041399}.

\bibitem{VandenFabio2006}
{\sc F.~A. Tal and E.~Vanden-Eijnden}, {\em Transition state theory and
  dynamical corrections in ergodic systems}, Nonlinearity, 19 (2006), p.~501,
  \href{http://dx.doi.org/10.1088/0951-7715/19/2/014}
  {doi:10.1088/0951-7715/19/2/014}.

\bibitem{TiaBur06}
{\sc T.~Tian and K.~Burrage}, {\em Stochastic models for regulatory networks of
  the genetic toggle switch.}, Proc.~Natl.~Acad.~Sci.~USA, 103 (2006),
  pp.~8372--8377, \href{http://dx.doi.org/10.1073/pnas.0507818103}
  {doi:10.1073/pnas.0507818103},
  \url{http://dx.doi.org/10.1073/pnas.0507818103}.

\bibitem{TouJerRut:2006}
{\sc A.~Touhami, M.~Jericho, and A.~D. Rutenberg}, {\em Temperature dependence
  of {MinD} oscillation in \emph{Escherichia coli}}, Mol.~Microb., 188 (2006),
  pp.~1279--1284, \href{http://dx.doi.org/10.1128/JB.00911-06}
  {doi:10.1128/JB.00911-06}.

\bibitem{Kam01}
{\sc N.~G. van Kampen}, {\em Stochastic Processes in Physics and Chemistry},
  Elsevier Science, 2001,
  \url{http://store.elsevier.com/Stochastic-Processes-in-Physics-and-Chemistry/N_G_-Van-Kampen/isbn-9780444529657/}.

\bibitem{WaAnDiCu}
{\sc J.~C. Walsh, C.~N. Angstmann, I.~G. Duggin, and P.~M.~G. Curmi}, {\em
  Molecular interactions of the {M}in protein system reproduce spatiotemporal
  patterning in growing and dividing \emph{Escherichia coli} cells}, PLoS ONE,
  10 (2015), p.~e0128148, \href{http://dx.doi.org/10.1371/journal.pone.0128148}
  {doi:10.1371/journal.pone.0128148}.

\bibitem{WeisseMiddletonHuisinga10}
{\sc A.~Y. Weisse, R.~H. Middleton, and W.~Huisinga}, {\em Quantifying
  uncertainty, variability and likelihood for ordinary differential equation
  models}, BMC Syst. Biol., 4 (2010), p.~144,
  \href{http://dx.doi.org/10.1186/1752-0509-4-144}
  {doi:10.1186/1752-0509-4-144}.

\bibitem{WinfreeGeometryOfBiologicaltime}
{\sc A.~T. Winfree}, {\em The Geometry of Biological Time}, Springer, 2001,
  \href{http://dx.doi.org/10.1007/978-1-4757-3484-3}
  {doi:10.1007/978-1-4757-3484-3}.

\end{thebibliography}

\end{document}